\documentclass{elsart}

\usepackage{amsmath}
\usepackage{amssymb}
\usepackage{amsfonts}
\usepackage{enumerate}
\usepackage{bibentry}

\begin{document}

\begin{frontmatter}



\title{A study of inverse trigonometric
integrals associated with three-variable Mahler measures, and some
related identities}


\author{Mathew D. Rogers}

\address{Department of Mathematics\\
         University of British Columbia\\
        Vancouver, BC, V6T-1Z2, Canada
        }

\begin{abstract}
We prove several identities relating three-variable Mahler
measures to integrals of inverse trigonometric functions.  After
deriving closed forms for most of these integrals, we obtain ten
explicit formulas for three-variable Mahler measures.  Several of
these results generalize formulas due to Condon and Lal\'in.  As a
corollary, we also obtain three $q$-series expansions for the
dilogarithm.
\end{abstract}

\end{frontmatter}
\def\theequation{\thesection.\arabic{equation}}
\def\init{\setcounter{equation}{0}}
\renewcommand\thesection{\arabic{section}}
\newtheorem{theorem}{Theorem}[section]
\newtheorem{proposition}[theorem]{Proposition}
\newtheorem{lemma}[theorem]{Lemma}
\newtheorem{definition}[theorem]{Definition}
\newtheorem{corollary}[theorem]{Corollary}
\newtheorem{conjecture}[theorem]{Conjecture}
\newtheorem{remark}{Remark}[section]
\def\rr{{\bf R}}
\def\cc{{\bf C}}
\def\nn{{\bf N}}
\def\zz{{\bf Z}}
\def\qq{{\bf Q}}
\def\uu{{\bf u}}
\def\cala{{\cal A}}
\def\calb{{\cal B}}
\def\calt{{\cal T}}
\def\jap#1{\langle #1\rangle}
\def\Li{{\operatorname{Li}}}
\def\m{{\operatorname{m}}}
\def\sn{{\operatorname{sn}}}
\def\cn{{\operatorname{cn}}}
\def\dn{{\operatorname{dn}}}
\def\am{{\operatorname{am}}}
\def\d{{\operatorname{d}}}
\def\K{{\operatorname{K}}}
\def\F{{\operatorname{F}}}
\def\T{{\operatorname{T}}}
\def\S{{\operatorname{S}}}
\def\TS{{\operatorname{TS}}}
\def\proof{{\it  Proof. }}
\def\half{\frac{1}{2}}
\def\eps{\epsilon}
\def\Re{\hbox{Re}\,}
\def\Im{\hbox{Im}\,}
\def\supp{\hbox{supp}\,}
\def\dist{\hbox{dist}\,}
\def\mod{\hbox{mod}\,}
\def\p{\partial}
\def\one{{\bf 1}}
\def\deg{\hbox{deg}\,}
\def\lcm{\hbox{lcm}\,}
\def\gcd{\hbox{gcd}\,}
\def\div{\,|\,}
\def\ndiv{\!\not|\,}
\def\aij{{\alpha_{ij}}}
\def\bij{{\beta_{ij}}}
%
\parindent=20pt
%
\section{Introduction}
\label{intro} \init
%
In this paper we will undertake a systematic study of each of the
inverse trigonometric integrals
\begin{equation*}
\begin{split}
\T(v,w)=\int_{0}^{1}\frac{\tan^{-1}(v x)\tan^{-1}(w x)}{x}\d x,\\
\S(v,w)=\int_{0}^{1}\frac{\sin^{-1}(v x)\sin^{-1}(w x)}{x}\d x,\\
\TS(v,w)=\int_{0}^{1}\frac{\tan^{-1}(v x)\sin^{-1}(w x)}{x}\d x.
\end{split}
\end{equation*}
This class of integrals arises when trying to find closed form
expressions for the Mahler measures of certain three-variable
polynomials.

Recall that the Mahler measure of an $n$-dimensional polynomial,
$P(x_1,\dots,x_n)$, can be defined by
\begin{equation*}
\m\left(P(x_1,\dots,x_n)\right)=\int_{0}^{1}\dots\int_{0}^{1}\log\big\vert
P\left(e^{2\pi i\theta_1},\dots,e^{2\pi i\theta_n}\right)\big\vert
\d\theta_1\dots\d\theta_n.
\end{equation*}
In the last few years, numerous papers have established explicit
formulas relating multi-variable Mahler measures to special
constants. Smyth proved the first result \cite{Bo0} with
\begin{equation*}
\m\left(1+x+y+z\right)=\frac{7}{2\pi^2}\zeta(3),
\end{equation*}
where the Riemann zeta function is defined by
$\zeta(s)=\sum_{n=1}^{\infty}\frac{1}{n^s}$.

In this paper, we will prove a number of new formulas relating
three-variable Mahler measures to the aforementioned trigonometric
integrals.  Many of our identities generalize previously known
results.  We will list a few of our main results in this
introductory section.

For our first example, we can use various properties of $\T(v,w)$
to show that
\begin{equation}\label{T(v,1/v) to mahler measure intro}
\begin{split}
\m&\left(1-
v^4\left(\frac{1-x}{1+x}\right)^2+\left(y+v^2\left(\frac{1-x}{1+x}\right)\right)^2
z\right)\\
&\qquad=\frac{4}{\pi}\int_{0}^{v}\frac{\tan^{-1}(u)}{u}\d
u-\frac{8}{\pi^2}\T\left(v,\frac{1}{v}\right)
+\frac{1}{2}\m\left(1-v^4\left(\frac{1-x}{1+x}\right)^2\right).
\end{split}
\end{equation}
This reduces to one of Lalin's formulas \cite{La1} when $v=1$:
\begin{equation}
\m\left((1+y)(1+z)+(1-z)(x-y)\right)=\frac{7}{2\pi^2}\zeta(3)+\frac{\log(2)}{2}.
\end{equation}

    We can use the double arcsine integral, $\S(v,w)$, to prove that
if $v\in[0,1]$:
\begin{equation}
\begin{split}
\m\left(v(1+x)+y+z\right)&=\frac{2}{\pi}\int_{0}^{v}\frac{\sin^{-1}(u)}{u}\d
u-\frac{4}{\pi^2}\S(v,1)\\
        &=\frac{4}{\pi^2}\left(\frac{\Li_3(v)-\Li_3(-v)}{2}\right).
\end{split}
\end{equation}
The second equality has been proved by Vandervelde \cite{Va}.
Slightly more complicated arguments lead to expressions that
include
\begin{equation}\label{intro S(1/2,1/2) formula}
\m\left(1-x^{1/6}+y+z\right)=\frac{2}{\pi}\int_{0}^{\frac{1}{2}}\frac{\sin^{-1}(u)}{u}\d
u-\frac{12}{\pi^2}\S\left(\frac{1}{2},\frac{1}{2}\right)
\end{equation}
This fractional Mahler measure is defined by
\begin{equation*}
\m\left(1-x^{1/6}+y+z\right)=\int_{0}^{1}\m\left(1-e^{2\pi i
u/6}+y+z\right)\d u,
\end{equation*}
notice that $\m\left(1-x^{1/6}+y+z\right)\not=\m(1-x+y+z)$. We can
simplify the right-hand side of Eq. \eqref{intro S(1/2,1/2)
formula} by either expressing
$\S\left(\frac{1}{2},\frac{1}{2}\right)$ as a linear combination
of $\operatorname{L}$-functions, or in terms of a famous binomial
sum:
\begin{equation*}
\S\left(\frac{1}{2},\frac{1}{2}\right)=\frac{1}{4}\sum_{n=1}^{\infty}\frac{1}{n^3{2n\choose
n}}.
\end{equation*}

    Condon \cite{Co} proved an identity that Boyd and Rodriguez
Villegas conjectured:
\begin{equation}
\m\left(1+x+(1-x)(y+z)\right)=\frac{28}{5\pi^2}\zeta(3).
\end{equation}
His proof also showed (in a slightly disguised form) that
\begin{equation}\label{condon intro}
\TS(2,1)=\frac{\pi}{2}\int_{0}^{2}\frac{\tan^{-1}(u)}{u}\d
u-\frac{7}{5}\zeta(3).
\end{equation}
We have generalized Condon's identity to show that
\begin{equation}
\m\left(1+x+\frac{v}{2}(1-x)(y+z)\right)=\frac{2}{\pi}\int_{0}^{v}\frac{\tan^{-1}(u)}{u}\d
u-\frac{4}{\pi^2}\TS(v,1),
\end{equation}
where Eq. \eqref{TS(v,1) closed form} expresses $\TS(v,1)$ in
terms of polylogarithms. We can use this result to prove a number
of new formulas, including:
\begin{equation}
\begin{split}
\m&\left(x+\frac{v^2}{4}(1+x)^2+\left(y+\frac{v}{2}(1+x)\right)^2
z\right)\\
&\qquad=\frac{2}{\pi}\int_{0}^{v}\frac{\tan^{-1}(u)}{u}\d
u-\frac{4}{\pi^2}\TS(v,1)
+\frac{1}{2}\m\left(x+\frac{v^2}{4}(1+x)^2\right).
\end{split}
\end{equation}
When $v=2$ this reduces to an interesting identity for $\zeta(3)$
and the golden ratio:
\begin{equation}
\m\left(x+(1+x)^2+(1+x+y)^2
z\right)=\frac{28}{5\pi^2}\zeta(3)+\log\left(\frac{1+\sqrt{5}}{2}\right).
\end{equation}

 We will show that all of the integrals $\TS(v,w)$, $\T(v,w)$, and $\S(v,w)$
have closed form expressions in terms of polylogarithms.  The
special case of $\TS(v,1)$ will warrant extra attention, as it is
related to an interesting family of binomial sums.  Our closed
forms are all derived through elementary methods.

\section{Preliminaries: A description of the method, and
some two dimensional Mahler measures}\label{mahler measure
integrals} \init

Although there are many conjectured formulas for multi-variable
Mahler measures, most are extremely difficult, if not impossible,
to prove.  Rather than attempting to prove any of these
conjectures, we will take an easier approach. By investigating
promising functions, and rewriting them as Mahler measures, we can
recover a number of useful formulas.

    Our first step was to determine a class of functions that we could relate
to Mahler's measure.  We chose the three integrals $\TS(v,w)$,
$\S(v,w)$, and $\T(v,w)$, based on Condon's evaluation of
$\TS(2,1)$, Eq. \eqref{condon intro}. Condon's formula naturally
suggested the existence of a generalized Mahler measure formula
involving $\TS(v,1)$. From there, it was a small step to consider
the similar functions $\TS(v,w)$, $\T(v,w)$, and $\S(v,w)$.

    We will use the following method to express $\TS(v,1)$,
$\S(v,1)$, and $\T(v,1/v)$ as three-variable Mahler measures.
First, a simple integration by parts changes each function into a
two-dimensional integral, containing either a nested arcsine or
arctangent integral. Recall that the following integrals define
the arctangent and arcsine integrals respectively:
\begin{align*}
&\int_{0}^{w}\frac{\tan^{-1}(u)}{u}\d u,
&&\int_{0}^{v}\frac{\sin^{-1}(u)}{u}\d u.
\end{align*}
A typical formula for $\TS(v,1)$, Eq. \eqref{TS(v,1) double
integral}, can be proved with little trouble:
\begin{equation*}
\TS(v,1)=\frac{\pi}{2}\int_{0}^{v}\frac{\tan^{-1}(u)}{u}\d
u-\int_{0}^{\pi/2}\int_{0}^{v\sin(\theta)}\frac{\tan^{-1}(z)}{z}\d
z\d\theta.
\end{equation*}
Next, substituting a two-dimensional Mahler measure for the nested
arctangent or arcsine integral will allow us to obtain a
three-dimensional Mahler measure evaluation.  Theorem \ref{condon
generalized theorem}, Proposition \ref{S related to mahler
proposition}, and Theorem \ref{T(v,1/v) mahler theorem} contain
our main results from using this method.

    Expressing the arcsine and arctangent
integrals in terms of Mahler's measure represents the main
difficulty in this approach.  In the remainder of this section we
will establish four two-variable Mahler measures for the
arctangent integral, and one two-variable Mahler measure for the
arcsine integral.

    Since many of our results involve polylogarithms, this will be
a good place to define the polylogarithm.
\begin{definition} If $|z|<1$, then the polylogarithm of order $k$ is defined by
\begin{equation*}
\Li_k(z)=\sum_{n=1}^{\infty}\frac{z^n}{n^k}.
\end{equation*}
We call $\Li_2(z)$ the dilogarithm, and we call $\Li_3(z)$ the
trilogarithm.
\end{definition}

    Theorem \ref{mahler pre integrals} requires a formula of Cassaigne and Maillot \cite{Ma}.
In particular, Cassaigne and Mallot showed that
\begin{equation*}
\pi\m(a+b x+c y)=\left\{\begin{array}{ll}
                        D\left(\frac{|a|}{|b|}e^{i \gamma}\right)+\alpha\log|a|+\beta\log|b|+\gamma\log|c|,&\text{ if ``$\triangle$"}\\
                        \pi\log\left(\max\left\{|a|,|b|,|c|\right\}\right), &\text{otherwise}.\\
                    \end{array}\right.\
\end{equation*}
The ``$\triangle$" condition states that $|a|$, $|b|$, and $|c|$
form the sides of a triangle.  If ``$\triangle$" is true, then
$\alpha$, $\beta$, and $\gamma$ denote the radian measures of the
angles opposite to the sides of length $|a|$, $|b|$, and $|c|$
respectively. In this formula, $D(z)$ denotes the Bloch-Wigner
dilogarithm. As usual,
\begin{equation*}
D(z)=\Im\left(\Li_2(z)\right)+\log|z|\arg(1-z).
\end{equation*}
Now that we have stated Cassaigne and Maillot's formula, we will
prove Theorem \ref{mahler pre integrals}.

\begin{theorem}\label{mahler pre integrals} If $0\le v\le 1$ and $w\ge 0$, then
\begin{align}
\int_{0}^{v}\frac{\sin^{-1}(u)}{u}\d u=&\frac{\pi}{2}\m(2v+y+z)\label{arcsine integral to mahler measure}\\
\int_{0}^{w}\frac{\tan^{-1}(u)}{u}\d
u=&\frac{\pi}{2}\m\left(1+w^2+(y+w)^2
z\right)-\frac{\pi}{4}\log\left(1+w^2\right)\label{arctangent
integral to mahler measure}
\end{align}
\end{theorem}
\begin{proof}To prove Eq. \eqref{arcsine integral to mahler
measure} first recall the usual formula for this arcsine integral,
\begin{equation}\label{arcsine integral formula}
\int_{0}^{v}\frac{\sin^{-1}(u)}{u}\d
u=\frac{1}{2}\Im\left(\Li_2\left(e^{2i\sin^{-1}(v)}\right)\right)+\sin^{-1}(v)\log(2v),
\end{equation}
which is valid whenever $0\le v\le 1$.

Now apply Cassaigne and Maillot's formula to $\m(2v+y+z)$; we are
in the ``$\triangle$" case since $0\le v\le 1$.  It follows from a
little trigonometry that
\begin{equation*}
\pi
\m(2v+y+z)=D\left(e^{2i\sin^{-1}(v)}\right)+2\sin^{-1}(v)\log(2v),
\end{equation*}
Since $\left|e^{2i\sin^{-1}(v)}\right|=1$,
$D\left(e^{2i\sin^{-1}(v)}\right)=\Im\left(\Li_2\left(e^{2i\sin^{-1}(v)}\right)\right)$,
hence we obtain
\begin{equation*}
\pi\m(2v+y+z)=\Im\left(\Li_2\left(e^{2i\sin^{-1}(v)}\right)\right)+2\sin^{-1}(v)\log(2v).
\end{equation*}
Comparing this last formula to Eq. \eqref{arcsine integral
formula}, we have
\begin{equation*}
\frac{\pi}{2}\m(2v+y+z)=\int_{0}^{v}\frac{\sin^{-1}(u)}{u}\d u.
\end{equation*}

To prove Eq. \eqref{arctangent integral to mahler measure} first
recall that if $0\le w\le 1$, then
\begin{equation*}
\int_{0}^{w}\frac{\tan^{-1}(u)}{u}\d u=\Im\left(\Li_2(i w)\right).
\end{equation*}
Next observe that by Cassaigne and Maillot's formula
\begin{align*}
\pi\m\left(\sqrt{1+w^2}+wy+z\right)&=D\left(e^{\pi
i/2}w\right)+\tan^{-1}(w)\log(w)+\frac{\pi}{2}\log\left(\sqrt{1+w^2}\right)\\
                &=\Im\left(\Li_2(iw)\right)+\frac{\pi}{4}\log\left(1+w^2\right).
\end{align*}
Making a change of variables in the Mahler measure, it is clear
that
\begin{align*}
\m\left(\sqrt{1+w^2}+wy+z\right)&=
        \frac{1}{2}\left\{\m\left(\sqrt{1+w^2}+(1+wy)iz\right)+\m\left(\sqrt{1+w^2}-(1+wy)iz\right)\right\}\\
        &=\frac{1}{2}\m\left(1+w^2+(1+wy)^2z^2\right)\\
        &=\frac{1}{2}\m\left(1+w^2+(y+w)^2 z\right).
\end{align*}
It follows that for $0\le w\le 1$ we have
\begin{equation*}
\int_{0}^{w}\frac{\tan^{-1}(u)}{u}\d
u=\frac{\pi}{2}\m\left(1+w^2+(y+w)^2z\right)-\frac{\pi}{4}\log\left(1+w^2\right).
\end{equation*}
We can extend this formula to the entire positive real line.
Suppose that $w=1/w'$ where $w'\ge 1$, then
\begin{align*}
\int_{0}^{1/w'}\frac{\tan^{-1}(u)}{u}\d
u&=\frac{\pi}{2}\m\left(1+\frac{1}{w'^2}+\left(y+\frac{1}{w'}\right)^2z\right)-\frac{\pi}{4}\log\left(1+\frac{1}{w'^2}\right)\\
&=\frac{\pi}{2}\m\left(1+w'^2+\left(y+w'\right)^2z\right)-\frac{\pi}{4}\log\left(1+w'^2\right)-\frac{\pi}{2}\log(w').
\end{align*}
Since the arctangent integral obeys the functional equation
\cite{Ra}
\begin{equation}\label{arctangent integral functional equation}
\int_{0}^{w'}\frac{\tan^{-1}(u)}{u}\d
u=\frac{\pi}{2}\log(w')+\int_{0}^{1/w'}\frac{\tan^{-1}(u)}{u}\d u,
\end{equation}
it follows that
\begin{align*}
\int_{0}^{w'}\frac{\tan^{-1}(u)}{u}\d u
=\frac{\pi}{2}\m\left(1+w'^2+\left(y+w'\right)^2z\right)-\frac{\pi}{4}\log\left(1+w'^2\right).
\end{align*}
Therefore Eq. \eqref{arctangent integral to mahler measure} holds
for all $w\ge 0$.$\blacksquare$
\end{proof}

    The next theorem proves that Eq. \eqref{arctangent
integral to mahler measure} is not unique. Using results from
Theorem \ref{integral table of jacobian integrals}, we can derive
three more Mahler measures for the arctangent integral.

\begin{theorem}\label{arctan extra mahlers theoreom} Suppose that $w\ge 0$, then
\begin{align}
 \int_{0}^{w}\frac{\tan^{-1}(u)}{u}\d
u=&\frac{\pi}{4}\m\left((1+w^2)(1+y)+w(1-y)(z+z^{-1})\right)\label{arctan extra mahler 1},\\
\int_{0}^{w}\frac{\tan^{-1}(u)}{u}\d
u=&\frac{\pi}{2}\m\left((y-y^{-1})+w(z+z^{-1})\right)\label{arctan
extra mahler 2},
\end{align}
\begin{equation}
\begin{split}
\int_{0}^{w}\frac{\tan^{-1}(u)}{u}\d
u=&\frac{\pi}{4}\m\left(\begin{split}&\left(4(1+y)^2-\left(z+z^{-1}\right)^2\right)(1+w^2)^2\\
                                     &+\left(z-z^{-1}\right)^2(1+y)^2(1-w^2)^2\end{split}\right)\\
  &-\frac{\pi}{4}\log(2)-\frac{\pi}{2}\log(1+w).\end{split}\label{arctan extra mahler 3}
\end{equation}
\end{theorem}
\begin{proof}  Since all three of these formulas have similar proofs, we will only prove
Eq. \eqref{arctan extra mahler 1} and Eq. \eqref{arctan extra
mahler 3}.  It is necessary to remark, that while Eq.
\eqref{arctan extra mahler 1} follows from Eq. \eqref{cn(u)/dn(u)
integral}, and Eq. \eqref{arctan extra mahler 3} follows from Eq.
\eqref{1/sn(u) integral}, we must start from Eq. \eqref{sn(u)
integral} to prove Eq. \eqref{arctan extra mahler 2}.

Now we will proceed with the proof of Eq. \eqref{arctan extra
mahler 1}.  From Eq. \eqref{cn(u)/dn(u) integral} we have
\begin{equation*}
\frac{\pi}{4k}\log\left(\frac{1+k}{1-k}\right)-\frac{2}{k}\Im\left(\Li_2(ir)\right)=\int_{0}^{1}\frac{\sin^{-1}(u)}{1-k^2
u^2}\d u,
\end{equation*}
where $k=\frac{2r}{1+r^2}$, and $0<k<1$.  After an integration by
parts this becomes
\begin{equation*}
\begin{split}
\frac{\pi}{4k}\log\left(\frac{1+k}{1-k}\right)-&\frac{2}{k}\Im\left(\Li_2(ir)\right)\\
&=\frac{\pi}{4k}\log\left(\frac{1+k}{1-k}\right)-
\frac{1}{2k}\int_{0}^{1}\log\left(\frac{1+k u}{1-k
u}\right)\frac{\d u}{\sqrt{1-u^2}}.
\end{split}
\end{equation*}
It follows immediately that
\begin{align*}
\Im\left(\Li_2(ir)\right)&=\frac{1}{4}\int_{0}^{\pi/2}\log\left(\frac{1+k
\sin(t)}{1-k \sin(t)}\right)\d t\\
&=\frac{1}{8}\int_{0}^{2\pi}\log^{+}\bigg{\vert}\frac{1+k\sin(t)}{1-k\sin(t)}\bigg{\vert}\d
t.
\end{align*}
Changing the ``$\log^{+}\vert \cdot\vert$" term into a Mahler
measure, which we can do by Jensen's formula, yields
\begin{align*}
\Im\left(\Li_2(i
r)\right)&=\frac{\pi}{4}\m\left(y+\frac{1+k\frac{z+z^{-1}}{2}}{1-k\frac{z+z^{-1}}{2}}\right).
\end{align*}
Since $k=\frac{2}{r+r^{-1}}$, we have
\begin{align*}
\Im\left(\Li_2(i
r)\right)=&\frac{\pi}{4}\m\left(y+\frac{r+r^{-1}+(z+z^{-1})}{r+r^{-1}-(z+z^{-1})}\right)\\
         =&\frac{\pi}{4}\m\left((1+y)(r+r^{-1})+(1-y)(z+z^{-1})\right)\\
           &-\frac{\pi}{4}\m\left(r+r^{-1}-(z+z^{-1})\right)\\
         =&\frac{\pi}{4}\m\left((1+y)(r+r^{-1})+(1-y)(z+z^{-1})\right)\\
          &-\frac{\pi}{4}\left(\log^{+}(r)+\log^{+}\left(\frac{1}{r}\right)\right)
\end{align*}
In order to substitute the arctangent integral for
$\Im\left(\Li_2(i r)\right)$, we will assume that $0<r<1$. With
this restriction, the formula becomes
\begin{align}\label{arctan 1st proof part formula}
\int_{0}^{r}\frac{\tan^{-1}(u)}{u}\d
u=&\frac{\pi}{4}\m\left((1+y)(r+r^{-1})+(1-y)(z+z^{-1})\right)\notag\\
    &-\frac{\pi}{4}\log\left(\frac{1}{r}\right)\notag\\
  =&\frac{\pi}{4}\m\left((1+y)(1+r^2)+r(1-y)(z+z^{-1})\right)
\end{align}
We can manually verify that Eq. \eqref{arctan 1st proof part
formula} holds when $r=0$ and $r=1$, and using Eq.
\eqref{arctangent integral functional equation} we can extend Eq.
\eqref{arctan 1st proof part formula} to all $r>1$. Therefore, Eq.
\eqref{arctan extra mahler 1} follows immediately.

    Next we will prove Eq. \eqref{arctan extra mahler 3}. Using Eq.
\eqref{sn(u) integral}, we can show that
\begin{equation*}
2\Im\left(\Li_2(i
p)\right)=\frac{\pi}{2}\log(p)+\int_{0}^{1}\frac{\sin^{-1}(u)}{u\sqrt{(1-u^2)(1-k^2
u^2)}}\d u,
\end{equation*}
where $k=\frac{1-p^2}{1+p^2}$, and $0<k<1$.  To satisfy this
restriction on $k$, we will assume that $0<p<1$.  After several
elementary simplifications, the right-hand side becomes
\begin{align*}
=&\frac{\pi}{2}\log(p)+\frac{\pi}{2}\log\left(1+\frac{1}{\sqrt{1-k^2}}\right)+
\int_{0}^{1}\log\left(1+\sqrt{\frac{1-u^2}{1-k^2
u^2}}\right)\frac{\d u}{\sqrt{(1-u^2)}}\\
=&\frac{\pi}{2}\log(p)+\frac{\pi}{2}\log\left(\frac{(1+p)^2}{2p}\right)+
\int_{0}^{\pi/2}\log\left(1+\frac{\cos(\theta)}{\sqrt{1-k^2
\sin^2(\theta)}}\right)\d\theta\\
=&\frac{\pi}{2}\log\left(\frac{(1+p)^2}{2}\right)+
\frac{1}{2}\int_{0}^{2\pi}\log^{+}\left\vert1+\frac{\cos(\theta)}{\sqrt{1-k^2
\sin^2(\theta)}}\right\vert\d\theta
\end{align*}
Since $\cos(\pi-\theta)=-\cos(\theta)$, we have
\begin{align*}
2\Im\left(\Li_2(i
p)\right)=&\frac{\pi}{2}\log\left(\frac{(1+p)^2}{2}\right)+\frac{1}{4}\int_{0}^{2\pi}\log^{+}\left\vert1+
\frac{\cos(\theta)}{\sqrt{1-k^2\sin^2(\theta)}}\right\vert\d\theta\\
&+\frac{1}{4}\int_{0}^{2\pi}\log^{+}\left\vert1-\frac{\cos(\theta)}{\sqrt{1-k^2
\sin^2(\theta)}}\right\vert\d\theta.
\end{align*}
Applying Jensen's formula yields
\begin{align*}
2\Im\left(\Li_2(i
p)\right)=&\frac{\pi}{2}\log\left(\frac{(1+p)^2}{2}\right)
+\frac{1}{4}\int_{0}^{2\pi}\m\left((1+y)^2-\frac{\cos^2(\theta)}{1-k^2
\sin^2(\theta)}\right)\d\theta\notag\\
        =&\frac{\pi}{2}\log\left(\frac{(1+p)^2}{2}\right)
            +\frac{\pi}{2}\m\left((1+y)^2-\frac{\left(z+z^{-1}\right)^2}{4+k^2\left(z-z^{-1}\right)^2}\right)\notag\\
        =&\frac{\pi}{2}\log\left(\frac{(1+p)^2}{2}\right)-\frac{\pi}{2}\m\left(4+k^2\left(z-z^{-1}\right)^2\right)\notag\\
            &+\frac{\pi}{2}\m\left(\left(4(1+y)^2-\left(z+z^{-1}\right)^2\right)+
                k^2(1+y)^2\left(z-z^{-1}\right)^2\right)
\end{align*}
We can simplify the one-dimensional Mahler measure as follows:
\begin{align*}
\m\left(4+k^2\left(z-z^{-1}\right)^2\right)
&=2\m\left(2+ik\left(z-z^{-1}\right)\right)\\
&=2\log\left(1+\sqrt{1-k^2}\right)\\
&=2\log\left(\frac{(1+p)^2}{1+p^2}\right).
\end{align*}
Eliminating $k$ yields
\begin{equation*}
\begin{split}
2\Im\left(\Li_2(i
p)\right)=&\frac{\pi}{2}\log\left(\frac{(1+p)^2}{2}\right)-\pi\log\left(\frac{(1+p)^2}{1+p^2}\right)\\
            &+\frac{\pi}{2}\m\left(\left(4(1+y)^2-\left(z+z^{-1}\right)^2\right)+
                \left(\frac{1-p^2}{1+p^2}\right)^2(1+y)^2\left(z-z^{-1}\right)^2\right)\\
         =&-\frac{\pi}{2}\log(2)-\pi\log(1+p)\\
            &+\frac{\pi}{2}\m\left(\begin{split}&\left(4(1+y)^2-\left(z+z^{-1}\right)^2\right)(1+p^2)^2\\
                                                &+(1+y)^2\left(z-z^{-1}\right)^2(1-p^2)^2\end{split}\right).
\end{split}
\end{equation*}
Since $0<p<1$, it follows that
\begin{equation}\label{arctan extra mahler 3 proof step}
\begin{split}
2\int_{0}^{p}\frac{\tan^{-1}(u)}{u}\d u=&\frac{\pi}{2}\m\left(\begin{split}&\left(4(1+y)^2-\left(z+z^{-1}\right)^2\right)(1+p^2)^2\\
                                                &+(1+y)^2\left(z-z^{-1}\right)^2(1-p^2)^2
                                    \end{split}\right)\\
                                    &-\frac{\pi}{2}\log(2)-\pi\log(1+p).
\end{split}
\end{equation}
It is relatively easy to verify that Eq. \eqref{arctan extra
mahler 3 proof step} holds when $p=0$ and $p=1$.  Using Eq.
\eqref{arctangent integral functional equation}, we can also
extend Eq. \eqref{arctan extra mahler 3 proof step} to $p>1$,
which completes the proof of Eq. \eqref{arctan extra mahler 3}.
 $\blacksquare$
\end{proof}
\section{Relations between $\TS(v,1)$ and Mahler's measure,
and a reduction of $\TS(v,w)$ to multiple polylogarithms} \label{A
closed form for TS section} \init
The first goal of this section is to establish five identities
relating $\TS(v,1)$ to three-variable Mahler measures. We will
prove these formulas in Theorem \ref{condon generalized theorem},
using the methods outlined in Section \ref{mahler measure
integrals}. Corollary \ref{condon mahlers reduced cor} examines a
few special cases of these results.

    Theorem \ref{TS(v,w) reduction theorem} accomplishes the second goal
of this section, which is to express $\TS(v,w)$ in terms of
multiple polylogarithms. This result, which appears to be new, is
stated in Eq. \eqref{TS(v,w) in terms of multi polylogs}.  The
importance of Eq. \eqref{TS(v,w) in terms of multi polylogs} lies
in its easy proof, and more importantly in the fact that it
immediately reduces $\TS(v,1)$ to multiple polylogarithms.
Finally, Proposition \ref{multi polylog reduction cor} will
demonstrate that the multiple polylogarithms in Eq. \eqref{TS(v,w)
in terms of multi polylogs} always reduce to standard
polylogarithms.

We will need the following simple lemma to prove Theorem \ref{ts
related to mahler}.
\begin{lemma}\label{ts related to mahler}
Assume that $v$ and $w$ are real numbers with $v>0$ and $w\in(0,1]$, then
\begin{align}
\begin{split}
\TS(v,w)=&\tan^{-1}(v)\int_{0}^{w}\frac{\sin^{-1}(z)}{z}\d z
-\int_{0}^{\tan^{-1}(v)}\int_{0}^{\frac{w}{v}\tan(\theta)}\frac{\sin^{-1}(z)}{z}\d
z\d \theta,\end{split}\label{ts with nested arcsine}\\
\begin{split}\TS(v,w)=&\sin^{-1}(w)\int_{0}^{v}\frac{\tan^{-1}(u)}{u}\d
u
-\int_{0}^{\sin^{-1}(w)}\int_{0}^{\frac{v}{w}\sin(\theta)}\frac{\tan^{-1}(z)}{z}\d
z\d\theta.\label{ts with nested arctan}
\end{split}
\end{align}
\end{lemma}
\begin{proof}  To prove Eq. \eqref{ts with nested arcsine} first integrate
$\TS(v,w)$ by parts to obtain:
\begin{equation*}
\begin{split}
 \TS(v,w)=&\tan^{-1}(v)\int_{0}^{w}\frac{\sin^{-1}(z)}{z}\d
z\\
&-\int_{0}^{1}\frac{\d}{\d u}\left(\tan^{-1}(v
u)\right)\int_{0}^{w u}\frac{\sin^{-1}(z)}{z}\d z\d u.
\end{split}
\end{equation*}
Making the $u$-substitution $\theta=\tan^{-1}(v u)$ we have:
\begin{equation*}
\TS(v,w)=\tan^{-1}(v)\int_{0}^{w}\frac{\sin^{-1}(z)}{z}\d
z-\int_{0}^{\tan^{-1}(v)}\int_{0}^{\frac{w}{v}\tan(\theta)}\frac{\sin^{-1}(z)}{z}\d
z\d \theta,
\end{equation*}
which completes the proof of the identity.

    The proof of Eq. \eqref{ts with nested arctan} follows in a
similar manner. $\blacksquare$
\end{proof}

The fact that Lemma \ref{ts related to mahler} expresses
$\TS(v,w)$ as a double integral in two different ways, makes
$\TS(v,w)$ more versatile than either $\S(v,w)$ or $\T(v,w)$.
These two different expansions will allow us to combine $\TS(v,w)$
with Mahler measures for both arctangent and arcsine integrals.

\begin{theorem}\label{condon generalized theorem} The following Mahler measures hold
whenever $v\ge 0$:{\allowdisplaybreaks
\begin{align}
\m&\left(1+x+\frac{v}{2}(1-x)(y+z)\right)=\frac{2}{\pi}\int_{0}^{v}\frac{\tan^{-1}(u)}{u}\d
u-\frac{4}{\pi^2}\TS(v,1)\label{condon generalized 1}\\
\begin{split}
\m&\left(x+\frac{v^2}{4}(1+x)^2+\left(y+\frac{v}{2}(1+x)\right)^2
z\right)\\
&\qquad=\frac{2}{\pi}\int_{0}^{v}\frac{\tan^{-1}(u)}{u}\d
u-\frac{4}{\pi^2}\TS(v,1)
+\frac{1}{2}\m\left(x+\frac{v^2}{4}(1+x)^2\right)
\end{split}\label{condon generalized
2}\\
\begin{split}
\m&\left((1+y)\left(1+\frac{v^2}{4}(x+x^{-1})^2\right)+
\frac{v}{2}(1-y)\left(x+x^{-1}\right)\left(z+z^{-1}\right)\right)\\
&\qquad=\frac{4}{\pi}\int_{0}^{v}\frac{\tan^{-1}(u)}{u}\d
u-\frac{8}{\pi^2}\TS(v,1)
\end{split}\label{condon generalized 3}\\
\begin{split}
\m&\left((z-z^{-1})+\frac{v}{2}(x+x^{-1})(y+y^{-1})\right)\\
&\qquad=\frac{2}{\pi}\int_{0}^{v}\frac{\tan^{-1}(u)}{u}\d
u-\frac{4}{\pi^2}\TS(v,1)
\end{split}\label{condon generalized 4}
\end{align}
\begin{equation}\label{condon generalized 5}
\begin{split}
\m&\left(\begin{split}&\left(4(1+y)^2-\left(z+z^{-1}\right)^2\right)\left(1+\frac{v^2}{4}\left(x+x^{-1}\right)^2\right)^2\\
&+\left(z-z^{-1}\right)^2(1+y)^2\left(1-\frac{v^2}{4}\left(x+x^{-1}\right)^2\right)^2\end{split}\right)\\
&\qquad=\frac{4}{\pi}\int_{0}^{v}\frac{\tan^{-1}(u)}{u}\d
u-\frac{8}{\pi^2} \TS(v,1)
+\frac{4}{\pi}\int_{0}^{\pi/2}\log\left(1+v\sin(\theta)\right)\d\theta\\
&\qquad\quad+\log(2)
\end{split}
\end{equation}}
\end{theorem}
\begin{proof}We will prove Eq. \eqref{condon generalized 1} first, since it has the most difficult
proof. Letting $w=1$ in Eq. \eqref{ts with nested arcsine} yields
\begin{equation*}
\TS(v,1)=\frac{\pi}{2}\log(2)\tan^{-1}(v)-\int_{0}^{\tan^{-1}(v)}
\int_{0}^{\tan(\theta)/v}\frac{\sin^{-1}(z)}{z}\d z\d \theta.
\end{equation*}
Since $0\le \frac{\tan(\theta)}{v}\le 1$, we may substitute Eq.
\eqref{arcsine integral to mahler measure} for the nested arcsine
integral to obtain
\begin{equation*}
\begin{split}
\TS(v,1)=&\frac{\pi}{2}\log(2)\tan^{-1}(v)-\frac{\pi}{2}\int_{0}^{\tan^{-1}(v)}\m\left(\frac{2}{v}\tan(\theta)+y+z\right)\d
\theta\\
=&\frac{\pi}{2}\log(2)\tan^{-1}(v)-\frac{\pi}{2}\int_{0}^{\pi/2}\m\left(\frac{2}{v}\tan(\theta)+y+z\right)\d
\theta\\
&+\frac{\pi}{2}\int_{\tan^{-1}(v)}^{\pi/2}\m\left(\frac{2}{v}\tan(\theta)+y+z\right)\d
\theta.
\end{split}
\end{equation*}
In the right-hand integral $\frac{\tan(\theta)}{v}\ge 1$, hence by
Cassaigne and Maillot's formula
\begin{equation*}
\m\left(\frac{2}{v}\tan(\theta)+y+z\right)=\log\left(\frac{2}{v}\tan(\theta)\right).
\end{equation*}
Substituting this result yields:
\begin{equation*}
\begin{split}
\TS(v,1)=&\frac{\pi}{2}\log(2)\tan^{-1}(v)+
\frac{\pi}{2}\int_{\tan^{-1}(v)}^{\pi/2}\log\left(\frac{2}{v}\tan(\theta)\right)\d\theta\\
&-\frac{\pi}{2}\int_{0}^{\pi/2}\m\left(\frac{2}{v}\tan(\theta)+y+z\right)\d
\theta\\
        =&\frac{\pi}{2}\int_{0}^{v}\frac{\tan^{-1}(u)}{u}\d
        u-\frac{\pi}{2}\int_{0}^{\pi/2}\m\left(\tan(\theta)+\frac{v}{2}(y+z)\right)\d\theta\\
        =&\frac{\pi}{2}\int_{0}^{v}\frac{\tan^{-1}(u)}{u}\d
        u-\frac{\pi^2}{4}\m\left(1+x+\frac{v}{2}(1-x)(y+z)\right).
\end{split}
\end{equation*}
Eq. \eqref{condon generalized 1} follows immediately from
rearranging this final identity.

    The proofs of equations \eqref{condon generalized 2} through \eqref{condon generalized 5} are virtually
identical, hence we will only prove Eq. \eqref{condon generalized
3}.  Letting $w=1$ in Eq. \eqref{ts with nested arctan}, we have
\begin{equation}\label{TS(v,1) double integral}
\TS(v,1)=\frac{\pi}{2}\int_{0}^{v}\frac{\tan^{-1}(u)}{u}\d
u-\int_{0}^{\pi/2}\int_{0}^{v\sin(\theta)}\frac{\tan^{-1}(z)}{z}\d
z\d\theta.
\end{equation}
Substituting Eq. \eqref{arctan extra mahler 1} for the nested
arctangent integral yields
\begin{align*}
\TS(v,1)=&\frac{\pi}{2}\int_{0}^{v}\frac{\tan^{-1}(u)}{u}\d
u\\
&-\frac{\pi}{4}\int_{0}^{\pi/2}\m\left((1+y)\left(1+v^2\sin^2(\theta)\right)+v\sin(\theta)(1-y)(z+z^{-1})\right)\d\theta\\
=&\frac{\pi}{2}\int_{0}^{v}\frac{\tan^{-1}(u)}{u}\d u\\
&-\frac{\pi^2}{8}\m\left((1+y)\left(1-\frac{v^2}{4}(x-x^{-1})^2\right)+\frac{v}{2i}(1-y)(x-x^{-1})(z+z^{-1})\right).
\end{align*}
Letting $x\rightarrow i x$, we obtain
\begin{align*}
\TS(v,1)=&\frac{\pi}{2}\int_{0}^{v}\frac{\tan^{-1}(u)}{u}\d u\\
&-\frac{\pi^2}{8}\m\left((1+y)\left(1+\frac{v^2}{4}(x+x^{-1})^2\right)+\frac{v}{2}(1-y)(x+x^{-1})(z+z^{-1})\right).
\end{align*}
Eq. \eqref{condon generalized 3} follows immediately from
rearranging this final equality.

    Finally, we will remark that the while Eq. \eqref{condon generalized 3}
follows from substituting Eq. \eqref{arctan extra mahler 1} into
Eq. \eqref{TS(v,1) double integral},  we must substitute Eq.
\eqref{arctangent integral to mahler measure} to prove Eq.
\eqref{condon generalized 2}, Eq. \eqref{condon generalized 4}
follows from substituting Eq. \eqref{arctan extra mahler 2}, and
Eq. \eqref{condon generalized 5} follows from substituting Eq.
\eqref{arctan extra mahler 3}. $\blacksquare$
\end{proof}
\begin{corollary}\label{condon mahlers reduced cor} The formulas in Theorem \ref{condon generalized
theorem} reduce, in order, to the following identities when $v=2$:
\begin{align}
\m&\left((1+x)+(1-x)(y+z)\right)=\frac{28}{5\pi^2}\zeta(3),\label{TS(2,1) to mahler 1}\\
\m&\left(x+(1+x)^2+(1+x+y)^2
z\right)=\frac{28}{5\pi^2}\zeta(3)+\log\left(\frac{1+\sqrt{5}}{2}\right),\label{TS(2,1) to mahler 2}\\
\begin{split}\m&\left(\left(1+x+z\right)\left(1+x^{-1}+z^{-1}\right)+
y(1+x-z)\left(1+x^{-1}-z^{-1}\right)\right)\\
&\qquad=\frac{56}{5\pi^2}\zeta(3),\end{split}\label{TS(2,1) to mahler 3}\\
\m&\left(\left(z-z^{-1}\right)+\left(x+x^{-1}\right)\left(y+y^{-1}\right)\right)=\frac{28}{5\pi^2}\zeta(3),\label{TS(2,1)
to mahler 4}
\end{align}
\begin{equation}\label{TS(2,1) to mahler 5}
\begin{split}
\m\left(\begin{split}&\left(4z(1+y)^2-(1+z)^2\right)\left(1+3x+x^2\right)^2\\
&+(1-z)^2(1+y)^2\left(1+x+x^2\right)^2\end{split}\right)
=&\frac{56}{5\pi^2}\zeta(3)+\frac{16}{3\pi}G \\
&+\log(2).
\end{split}
\end{equation}
In Eq. \eqref{TS(2,1) to mahler 5}, and throughout the rest of the
paper, $G$ denotes Catalan's constant. In particular,
$G=1-\frac{1}{3^2}+\frac{1}{5^2}-\frac{1}{7^2}\dots$
\end{corollary}
\begin{proof}  As we have already stated, Condon proved Eq.
\eqref{TS(2,1) to mahler 1} in \cite{Co}.  His proof also showed
that \[\TS(2,1)=\frac{\pi}{2}\int_{0}^{2}\frac{\tan^{-1}(u)}{u}\d
u-\frac{7}{5}\zeta(3).\] Using this formula, equations
\eqref{TS(2,1) to mahler 2} through \eqref{TS(2,1) to mahler 5}
follow immediately from Theorem \ref{condon generalized
theorem}.$\blacksquare$
\end{proof}

Theorem \ref{condon generalized theorem} shows that we can obtain
closed forms for several three-variable Mahler measures by
reducing $\TS(v,1)$ to polylogarithms.  We have proved a
convenient closed form for $\TS(v,1)$ in Eq. \eqref{TS(v,1) closed
form}.
 Corollary \ref{TS(v,1) and TS(2,1) cor} also shows that this
closed form immediately implies Condon's evaluation of $\TS(2,1)$.
We will postpone further discussion of Eq. \eqref{TS(v,1) closed
form} until Section \ref{TS(v,1) evaluation section}.

    We will devote the remainder of this section to deriving a closed
form for $\TS(v,w)$ in terms of multiple polylogarithms.  For
convenience, we will use a slightly non-standard notation for our
multiple polylogarithms.

\begin{definition}Define $\F_j(x)$ by
\begin{equation*}
\F_j(x)=\sum_{n=0}^{\infty}\frac{x^{2n+1}}{(2n+1)^j}=\frac{\Li_j(x)-\Li_j(-x)}{2},
\end{equation*}
and define $\F_{j,k}(x,y)$ by
\begin{equation*}
\F_{j,k}(x,y)=\sum_{n=0}^{\infty}\frac{x^{2n+1}}{(2n+1)^{j}}\sum_{m=0}^{n}\frac{y^{2m+1}}{(2m+1)^{k}}.
\end{equation*}
\end{definition}

We will employ this notation throughout the rest of the paper.

\begin{theorem}\label{TS(v,w) reduction theorem}If $\frac{v}{w}\not\in
(-i\infty,-i]\cup[i,i\infty)$ and $w\in [-1,1]$, then we can
express $\TS(v,w)$ in terms of multiple polylogarithms. Let
$R=\frac{\frac{v}{w}}{1+\sqrt{1+\left(\frac{v}{w}\right)^2}}$, and
let $S=iw+\sqrt{1-w^2}$, then
\begin{equation}\label{TS(v,w) in terms of multi polylogs}
\begin{split}
\TS(v,w)=&2\F_3(R)-\F_3(R S)-\F_3(R/S)-4\F_{1,2}(R,1)\\
        &+2\F_{1,2}(R,S)+2\F_{1,2}(R,1/S)\\
        &+i\sin^{-1}(w)\left\{\F_2(R
        S)-\F_2(R/S)\right.\\
        &\qquad\qquad\qquad\left.-2\F_{1,1}(R,S)+2\F_{1,1}(R,1/S)\right\}.
\end{split}
\end{equation}
\end{theorem}
\begin{proof}  First note that by $u$-substitution
\begin{equation}\label{ts after u-substitution}
\TS(v,w)=\int_{0}^{\sin^{-1}(w)}\tan^{-1}\left(\frac{v}{w}\sin(\theta)\right)\cot(\theta)\theta
\d \theta.
\end{equation}
Since $w\in [-1,1]$, it follows that our path of integration is
along the real axis.  Next substitute the Fourier series
\begin{equation}\label{arctan fourier series}
\tan^{-1}\left(\frac{v}{w}\sin(\theta)\right)=2\sum_{n=0}^{\infty}\frac{R^{2n+1}}{2n+1}\sin\left((2n+1)\theta\right),
\end{equation}
into Eq. \eqref{ts after u-substitution}.  Swapping the order of
summation and integration, we have
\begin{equation*}
\TS(v,w)=2\sum_{n=0}^{\infty}\frac{R^{2n+1}}{2n+1}\int_{0}^{\sin^{-1}(w)}\sin\left((2n+1)\theta\right)\cot(\theta)\theta
\d \theta.
\end{equation*}
Uniform convergence justifies this interchange of summation and
integration.  In particular, Eq. \eqref{arctan fourier series}
converges uniformly whenever $|R|<1$ and $\theta \in \mathbb{R}$.
It is easy to show that $|R|<1$ except when $\frac{v}{w}\in
(-i\infty,-i]\cup[i,i\infty)$, in which case $|R|=1$.  If $|R|=1$,
then Eq. \eqref{arctan fourier series} no longer converges
uniformly, and hence the following arguments do not apply.

 Evaluating the nested integral
yields
\begin{equation}\label{TS almost done}
\begin{split}
\TS(v,w)=4\sum_{n=0}^{\infty}\frac{R^{2n+1}}{2n+1}\bigg\{&\sin^{-1}(w){\sum_{k=0}^{n}}^{\prime}\frac{\sin\left((2k+1)\sin^{-1}(w)\right)}{2k+1}\\
&-{\sum_{k=0}^{n}}^{\prime}\frac{1-\cos\left((2k+1)\sin^{-1}(w)\right)}{(2k+1)^2}\bigg\},
\end{split}
\end{equation}
where ${\sum\limits_{k=0}^{n}}^\prime
a_k=a_0+\dots+a_{n-1}+\frac{a_n}{2}$.  Simplifying Eq. \eqref{TS
almost done} completes our proof. $\blacksquare$
\end{proof}

    Eq. \eqref{TS(v,w) in terms of multi polylogs} deserves a
few remarks, since it is a fairly general result.  Firstly,
observe that a closer analysis of Eq. \eqref{arctan fourier
series} would probably allow us to relax the restriction that
$w\in[-1,1]$. Secondly, Eq. \eqref{TS(v,w) in terms of multi
polylogs} most likely has applications beyond the scope of this
paper.  For example, we can use Eq. \eqref{TS(v,w) in terms of
multi polylogs} to reduce the right-hand side of the following
equation
\begin{equation}
\begin{split}
&\sum_{n=1}^{\infty}\frac{(-1)^n}{(2n+1)^2}{2n\choose
n}\left(\frac{w}{2}\right)^{2n+1}\sum_{k=1}^{n}\frac{(-1)^{k+1}}{k}\\
&\qquad=\TS(1,w)-\frac{\pi}{4}\int_{0}^{w}\frac{\sin^{-1}(t)}{t}\d
t+\frac{\log(2)}{2}\int_{0}^{w}\frac{\sinh^{-1}(t)}{t}\d t,
\end{split}
\end{equation}
to multiple polylogarithms.

    We can use the final result of this section, Proposition
\ref{multi polylog reduction cor}, to reduce $\TS(v,w)$ to regular
polylogarithms.  This proposition allows us to equate $\TS(v,w)$
with a formula involving around twenty trilogarithms. While a
clever usage of trilogarithmic functional equations might simplify
this result, it seems more convenient to simply leave Eq.
\eqref{TS(v,w) in terms of multi polylogs} in its current form.

\begin{proposition}\label{multi polylog reduction cor} The functions $\F_{1,1}(x,y)$ and $\F_{1,2}(x,y)$ can be expressed
in terms of polylogarithms,  we have:
\begin{equation}\label{F11(x,y) closed form}
\begin{split}
4\F_{1,1}(x,y)=&\Li_2\left(\frac{x(1+y)}{1+x}\right)-\Li_2\left(\frac{x(1-y)}{1+x}\right)\\
&-\Li_2\left(\frac{-x(1+y)}{1-x}\right)+\Li_2\left(\frac{-x(1-y)}{1-x}\right).
\end{split}
\end{equation}
To reduce $\F_{1,2}(x,y)$ to polylogarithms, apply Lewin's
formula, Eq. \eqref{Lewin's Formula}, four times to the following
identity:
\begin{equation}\label{F12(x,y) integral}
\begin{split}
\F_{1,2}(x,y)=&\F_{3}(x
y)-\frac{1}{2}\log\left(1-x^2\right)\F_{2}(x
y)\\
&+\frac{1}{4}\int_{0}^{x}\frac{\log\left(1-u^2\right)\log\left(\frac{1+y
u}{1-y u}\right)}{u}\d u.
\end{split}
\end{equation}
\end{proposition}
\begin{proof}To prove Eq. \eqref{F12(x,y) integral}, first
swap the order of summation to obtain
\begin{equation*}
\F_{1,2}(x,y)=\F_{3}(x
y)+\F_{1}(x)\F_{2}(y)-\sum_{n=0}^{\infty}\frac{y^{2n+1}}{(2n+1)^2}\sum_{k=0}^{n}\frac{x^{2k+1}}{2k+1}.
\end{equation*}
Substituting an integral for the nested sum yields
\begin{align*}
\F_{1,2}(x,y)&=\F_{3}(x
y)+\F_{1}(x)\F_{2}(y)-\sum_{n=0}^{\infty}\frac{y^{2n+1}}{(2n+1)^2}\int_{0}^{x}\frac{1-u^{2n+2}}{1-u^2}\d
u\\
&=\F_{3}(x y)+\int_{0}^{x}\frac{u}{1-u^2}\F_{2}(y u)\d u.
\end{align*}
Integrating by parts, the identity becomes
\begin{equation*}
\begin{split}
\F_{1,2}(x,y)=&\F_{3}(x y)-\frac{1}{2}\log(1-x^2)\F_{2}(x
y)\\
&+\frac{1}{4}\int_{0}^{x}\frac{\log\left(1-u^2\right)\log\left(\frac{1+
y u}{1-y u}\right)}{u}\d u,
\end{split}
\end{equation*}
which completes the proof of Eq. \eqref{F12(x,y) integral}.

    We can verify Eq. \eqref{F11(x,y) closed form} by
differentiating each side of the equation with respect to $y$.
 $\blacksquare$
\end{proof}

Finally, observe that we can obtain simple closed forms for
$\F_{1,2}(x,1)$ and $\F_{2,1}(1,x)$ from Eq. \eqref{full
evaluation of lambda21(x)}.

\section{An evaluation of $\TS(v,1)$ using infinite series}
\label{TS(v,1) evaluation section} \init

    This evaluation of $\TS(v,1)$ generalizes a theorem due
to Condon. Condon proved a formula that Boyd and Rodriguez
Villegas conjectured:
\begin{equation*}
m\big(1+x+(1-x)(y+z)\big)=\frac{28}{5\pi^2}\zeta(3).
\end{equation*}
Condon's result is equivalent to evaluating $\TS(2,1)$ in closed
form. As Theorem \ref{condon generalized theorem} has shown,
generalizing this Mahler measure depends on finding a closed form
for $\TS(v,1)$. Eq. \eqref{TS(v,1) closed form} accomplishes this
goal by expressing $\TS(v,1)$ in terms of polylogarithms.

    This calculation of $\TS(v,1)$ is based on several series
transformations.  The first step is to expand $\TS(v,1)$ in a
Taylor series; observe that the following formula holds whenever
$|v|<1$:
\begin{equation}\label{taylor expansion}
\TS(v,1)=\frac{\pi}{2}\sum_{k=0}^{\infty}\frac{(-1)^k}{(2k+1)^2}
v^{2k+1}-\frac{1}{2}\sum_{k=0}^{\infty}\frac{(-1)^k}{(2k+1)^3}\frac{(2v)^{2k+1}}{{2k\choose
k}}.
\end{equation}
We can easily prove Eq. \eqref{taylor expansion} by starting from
Eq. \eqref{TS(v,1) double integral}.  Formula \eqref{taylor
expansion} shows that $\TS(v,1)$ is analytic in the open unit
disk. Unfortunately Eq. \eqref{taylor expansion} does not converge
when $v=2$, and hence it can not be used to calculate $\TS(2,1)$.
It will be necessary to find an analytic continuation of
$\TS(v,1)$ in order to carry out any useful computations.

The following family of functions will play a crucial role in our
calculations.
\begin{definition}
Define $h_{n}(v)$ by the infinite series,
\begin{equation}\label{define h n(v)}
h_{n}(v)=\sum_{k=0}^{\infty}\frac{(-1)^k}{(2k+1)^n}\frac{(2v)^{2k+1}}{{2k\choose
k}}.
\end{equation}
\end{definition}

Using the definition of $h_3(v)$, combined with the identity
\begin{equation*}
\sum_{k=0}^{\infty}\frac{(-1)^k}{(2k+1)^2}v^{2k+1}=\int_{0}^{v}\frac{\tan^{-1}(u)}{u}du,
\end{equation*}
it follows that Eq. \eqref{taylor expansion} can be rewritten as
\begin{equation}\label{rewritten taylor series}
\TS(v,1)=\frac{\pi}{2}\int_{0}^{v}\frac{\tan^{-1}(u)}{u}du-\frac{1}{2}h_3(v).
\end{equation}

    Finding a closed form for $\TS(v,1)$ we will entail
finding a closed form for $h_3(v)$.  Theorem \ref{evaluate h3}
accomplishes this goal, however first we need to prove several
auxiliary lemmas.  The idea behind our proof is very simple: first
we will find a closed form for $h_2(v)$ and then integrate it to
find a closed form for $h_3(v)$.

    Batir recently used this method in an interesting paper \cite{Ba} to
obtain a formula that is equivalent to Eq. \eqref{batir equation}.
Unfortunately Batir seems to have missed Eq. \eqref{h3 evaluated},
so we will provide a full derivation of this important result.

\begin{lemma}\label{h 2 identity lemma} The function $h_2(v)$ is analytic if $v\not \in
(-i\infty,-i]\cup[i,i\infty)$. Furthermore, we can express
$h_2(v)$ in terms of the dilogarithm,
\begin{equation}\label{h 2 identity}
\begin{split}
h_2(v)&=4\sum_{k=0}^{\infty}\frac{1}{(2k+1)^2}\bigg(\frac{v}{1+\sqrt{1+v^2}}\bigg)^{2k+1}\\
      &=2\Li_2\bigg(\frac{v}{1+\sqrt{1+v^2}}\bigg)-2\Li_2\bigg(\frac{-v}{1+\sqrt{1+v^2}}\bigg).
\end{split}
\end{equation}
\end{lemma}
\begin{proof}
We use the following elementary identity to prove Eq. \eqref{h 2
identity},
\begin{equation}\label{finite binomial sum}
\frac{2^{4k}}{(2k+1)^2{2k\choose
k}}=\sum_{j=0}^{\infty}\frac{(-1)^j}{2j+1}\frac{(2k)!}{(k+j+1)!(k-j)!}.
\end{equation}
Substituting Eq. \eqref{finite binomial sum} into the definition
of $h_2(v)$, we have
\begin{equation*}
\begin{split}
h_2(v)&=\sum_{k=0}^{\infty}\frac{(-1)^k}{(2k+1)^2}\frac{(2v)^{2k+1}}{{2k\choose
k}}\\
      &=4\sum_{k=0}^{\infty}(-1)^k\big(\frac{v}{2}\big)^{2k+1}\sum_{j=0}^{\infty}
      \frac{(-1)^j}{2j+1}\frac{(2k)!}{(k+j+1)!(k-j)!}.
\end{split}
\end{equation*}
If we assume that $|v|<1$, then the series converges uniformly,
hence we may swap the order of summation to obtain
\begin{equation*}
\begin{split}
h_2(v)&=4\sum_{j=0}^{\infty}\frac{1}{2j+1}\sum_{k=0}^{\infty}(-1)^{k+2j}
\frac{(2k+2j)!}{(k+2j+1)!k!}\left(\frac{v}{2}\right)^{2k+2j+1}\\
    &=4\sum_{j=0}^{\infty}\frac{1}{(2j+1)^2}\left(\frac{v}{2}\right)^{2j+1}
        \sum_{k=0}^{\infty}\frac{(j+\frac{1}{2})_k(j+1)_k}{(2j+2)_k}\frac{(-v^2)^k}{k!},
\end{split}
\end{equation*}
where $(x)_n=\frac{\Gamma(x+n)}{\Gamma(x)}$.  But then we have
\begin{equation*}
\begin{split}
h_2(v)&=4\sum_{j=0}^{\infty}\frac{1}{(2j+1)^2}\left(\frac{v}{2}\right)^{2j+1}
        {_2F_1}\left[\substack{j+\frac{1}{2},j+1\\2j+2}\big|-v^2\right],
\end{split}
\end{equation*}
where $_2F_1\left[\substack{a,b\\c}\big|x\right]$ is the usual
hypergeometric function.  A standard hypergeometric identity
\cite{Gr} shows that
\begin{equation*}
{_2F_1}\big[\substack{j+\frac{1}{2},j+1\\2j+2}\big|-v^2\big]
=\frac{2^{2j+1}}{(1+\sqrt{1+v^2})^{2j+1}},
\end{equation*}
from which we obtain
\begin{equation*}
h_2(v)=4\sum_{j=0}^{\infty}\frac{1}{(2j+1)^2}\bigg(\frac{v}{1+\sqrt{1+v^2}}\bigg)^{2j+1},
\end{equation*}
concluding the proof of the identity.

    We can use Eq. \eqref{h 2 identity} to analytically continue
$h_2(v)$ to a larger domain.  Recall that $\Li_2(r)-\Li_2(-r)$ is
analytic whenever $r\not\in (-\infty,-1]\cup[1,\infty)$, and
$\frac{v}{1+\sqrt{1+v^2}}$ is analytic whenever $v \not \in
(-i\infty,-i]\cup[i,i\infty)$. Since we have already assumed that
$v \not \in (-i\infty,-i]\cup[i,i\infty)$, we simply have to show
that the range of $r=\frac{v}{1+\sqrt{1+v^2}}$ does not intersect
the set $\left\{(-\infty,-1]\cup [1,\infty)\right\}$.

    Some elementary calculus shows that $|r|=\left|\frac{v}{1+\sqrt{1+v^2}}\right|\le
1$ for all $v\in \mathbb{C}$, with equality occurring only when
$v\in(-i\infty,-i]\cup[i,i\infty)$.  It follows that $h_2(v)$ is
analytic on
$\mathbb{C}-\left\{(-i\infty,-i]\cup[i,i\infty)\right\}$.$\blacksquare$
\end{proof}

    Since we have now expressed $h_2(v)$ in terms of dilogarithms, we can
find a closed form for $h_1(v)$ by differentiating Eq. \eqref{h 2
identity}:
\begin{equation}
h_1(v)=\frac{2}{\sqrt{1+v^2}}\log\big(v+\sqrt{1+v^2}\big).
\end{equation}

In Theorem \ref{evaluate h3}, we will integrate Eq. \eqref{h 2
identity} to find a closed form for $h_3(v)$ involving
trilogarithms.  To prove this theorem, we first need to establish
two lemmas. Lemma \ref{integral for roots} evaluates a necessary
integral, while Lemma \ref{double polylog evaluation} expresses
$\F_{2,1}(1,x)$ in terms of polylogarithms.

\begin{lemma}\label{integral for roots} If $j\ge 0$ is an integer, and $r=\frac{v}{1+\sqrt{1+v^2}}$, then we
have the following identity:
\begin{equation}
\int_{0}^{v}\frac{1}{u}\bigg(\frac{u}{1+\sqrt{1+u^2}}\bigg)^{2j+1}\d
u =\log\left(\frac{1+r}{1-r}\right)+\frac{r^{2j+1}}{2j+1}
-2\sum_{k=0}^{j}\frac{r^{2k+1}}{2k+1}.
\end{equation}
\end{lemma}
\begin{proof} To evaluate the integral
\begin{equation*}
w_j(v)=\int_{0}^{v}\frac{1}{u}\bigg(\frac{u}{1+\sqrt{1+u^2}}\bigg)^{2j+1}\d
u,
\end{equation*}
first make the substitution $z=\frac{u}{1+\sqrt{1+u^2}}$. In
particular we can show that $u=\frac{2z}{1-z^2}$ and $\frac{\d
u}{\d z}=2\frac{(1+z^2)}{(1-z^2)^2}$. Therefore we have
\begin{equation*}
\begin{split}
w_j(v)&=\int_{0}^{r}z^{2j}\left(\frac{1+z^2}{1-z^2}\right)\d z\\
      &=\int_{0}^{r}\frac{2}{1-z^2}\d z
      -\int_{0}^{r}\frac{1-z^{2j}}{1-z^2}\d z
      -\int_{0}^{r}\frac{1-z^{2j+2}}{1-z^2}\d z.
\end{split}
\end{equation*}
Next substitute the geometric series
$\frac{1-z^{2j}}{1-z^2}=\sum_{k=0}^{j-1}z^{2k}$ into each of the
right-hand integrals, and swap the order of summation and
integration to obtain
\begin{equation*}
\begin{split}
w_j(v)=&\int_{0}^{r}\frac{2}{1-z^2}\d z
       -\sum_{k=0}^{j-1}\frac{r^{2k+1}}{2k+1}-\sum_{k=0}^{j}\frac{r^{2k+1}}{2k+1}\\
      =&\log\left(\frac{1+r}{1-r}\right)+\frac{r^{2j+1}}{2j+1}
      -2\sum_{k=0}^{j}\frac{r^{2k+1}}{2k+1}.
\end{split}
\end{equation*}$\blacksquare$
\end{proof}

\begin{lemma}\label{double polylog evaluation} The following double
polylogarithm
\begin{equation}\label{double polylog}
\F_{2,1}(1,x)=\sum_{n=0}^{\infty}\frac{1}{(2n+1)^2}\sum_{k=0}^{n}\frac{x^{2k+1}}{2k+1}
\end{equation}
can be evaluated in closed form.  If $|x|< 1$,
\begin{equation}\label{full evaluation of lambda21(x)}
\begin{split}
8\F_{2,1}(1,x)=&4\Li_3(x)-\Li_3(x^2)-4\Li_3(1-x)-4\Li_3\left(\frac{x}{1+x}\right)+4\zeta(3)\\
                    &+\log\left(\frac{1+x}{1-x}\right)\Li_2(x^2)+\frac{\pi^2}{2}\log(1+x)+\frac{\pi^2}{6}\log(1-x)\\
                    &+\frac{2}{3}\log^3(1+x)-2\log(x)\log^2(1-x)
\end{split}
\end{equation}
\end{lemma}
\begin{proof}We will verify Eq. \eqref{full evaluation of lambda21(x)}
by differentiating each side of the identity.  First observe that
the infinite series in Eq. \eqref{double polylog} converges
uniformly whenever $|x|\le 1$, hence term by term differentiation
is justified at all points in the open unit disk. It follows that
\begin{equation}\label{lambda21(x) derivative}
\begin{split}
\frac{\d}{\d
x}\F_{2,1}(1,x)=&\sum_{n=0}^{\infty}\frac{1}{(2n+1)^2}\left(\frac{1-x^{2n+2}}{1-x^2}\right)\\
                    =&\frac{\pi^2}{8}\left(\frac{1}{1-x^2}\right)-\frac{x}{1-x^2}\left(\Li_2(x)-\frac{1}{4}\Li_2(x^2)\right),
\end{split}
\end{equation}
whenever $|x|<1$.

    Let $\varphi(x)$ denote the right-hand side of
Eq. \eqref{full evaluation of lambda21(x)}.  Taking the derivative
of $\varphi(x)$ we obtain:
\begin{equation}\label{double polylog varphi derivative}
\begin{split}
\frac{\d\varphi}{\d
x}=&\frac{4}{x}\Li_2(x)-\frac{2}{x}\Li_2(x^2)+\frac{4}{1-x}\Li_2(1-x)\\
    &-4\left(\frac{1}{x}-\frac{1}{1+x}\right)\Li_2\left(\frac{x}{1+x}\right)
    +\frac{2}{1-x^2}\Li_2(x^2)\\
   &-\frac{2}{x}\left(\log^2(1+x)-\log^2(1-x)\right)+\frac{\pi^2}{2}\left(\frac{1}{1+x}\right)
   -\frac{\pi^2}{6}\left(\frac{1}{1-x}\right)\\
   &+\frac{2}{1+x}\log^2(1+x)-\frac{2}{x}\log^2(1-x)+\frac{4}{1-x}\log(x)\log(1-x)
\end{split}
\end{equation}
We can simplify Eq. \eqref{double polylog varphi derivative} by
eliminating $\Li_2(1-x)$ and $\Li_2\left(\frac{x}{1+x}\right)$
with the functional equations:
\begin{align*}
\Li_2(1-x)&=\frac{\pi^2}{6}-\log(x)\log(1-x)-\Li_2\left(x\right),\\
\Li_2\left(\frac{x}{1+x}\right)&=-\frac{1}{2}\log^2(1+x)+\Li_2(x)-\frac{1}{2}\Li_2(x^2).
\end{align*}
Substituting these identities into Eq. \eqref{double polylog
varphi derivative} and simplifying, we are left with
\begin{align*}
\frac{\d\varphi}{\d
x}&=\pi^2\left(\frac{1}{1-x^2}\right)-\frac{8x}{1-x^2}\left(\Li_2(x)-\frac{1}{4}\Li_2(x^2)\right)\\
  &=\frac{\d}{\d x}\left\{8\F_{2,1}(1,x)\right\}.
\end{align*}
Eq. \eqref{lambda21(x) derivative} justifies this final step.
Since the derivatives of $8\F_{2,1}(1,x)$ and $\varphi(x)$ are
equal on the open unit disk, and since both functions vanish at
zero, we may conclude that
$8\F_{2,1}(1,x)=\varphi(x)$.$\blacksquare$
\end{proof}

  The proof of Eq. \eqref{full evaluation of lambda21(x)} requires a remark.
Despite the fact that the right-hand side of Eq. \eqref{full
evaluation of lambda21(x)} is single valued and analytic whenever
$|x|<1$, the individual terms involving $\Li_3(1-x)$ and $\log(x)$
are multivalued for $x\in (-1,0)$.  To avoid all ambiguity, we can
simply use $\F_{2,1}(1,x)=-\F_{2,1}(1,-x)$ to calculate the
function at negative real arguments .

\begin{theorem}\label{evaluate h3}
The function $h_3(v)$ is analytic on
$\mathbb{C}-\left\{(-i\infty,-i]\cup[i,i\infty)\right\}$. If
$v\not\in (-i\infty,-i]\cup[i,i\infty)$, then $h_3(v)$ can be
expressed in terms of polylogarithms.  Let
$r=\frac{v}{1+\sqrt{1+v^2}}$, then
\begin{equation}\label{h3 evaluated}
\begin{split}
h_3(v)=&\frac{1}{2}\Li_3(r^2)+4\Li_3(1-r)+4\Li_3\left(\frac{r}{1+r}\right)-4\zeta(3)\\
       &-\log\left(\frac{1+r}{1-r}\right)\Li_2(r^2)
       -\frac{2\pi^2}{3}\log(1-r)-\frac{2}{3}\log^3(1+r)\\
       &+2\log(r)\log^2(1-r).
\end{split}
\end{equation}
We can recover an equivalent form of Condon's identity by letting
$v=2$:
\begin{equation}\label{h_3(2) evaluated}
h_3(2)=\frac{14}{5}\zeta(3).
\end{equation}
\end{theorem}
\begin{proof} This proof is very simple since we have already
completed all of the hard computations.  Observe from Eq.
\eqref{define h n(v)} that if $|v|<1$,
\begin{equation}\label{h3 from h2}
h_3(v)=\int_{0}^{v}\frac{h_2(u)}{u}\d u.
\end{equation}
Lemma \ref{h 2 identity lemma} shows that $h_2(v)$ is analytic
provided that $v\not\in (-i\infty,-i]\cup[i,i\infty)$. If we
assume that the path of integration does not pass through either
of these branch cuts, then it is easy to see that Eq. \eqref{h3
from h2} provides an analytic continuation of $h_3(v)$ to
$\mathbb{C}-\left\{ (-i\infty,-i]\cup[i,i\infty)\right\}$.

    Next we will prove Eq. \eqref{h3 evaluated}.  Substituting Eq. \eqref{h 2 identity}
into Eq. \eqref{h3 from h2} yields an infinite series for $h_3(v)$
that is valid whenever $v\not\in(-i\infty,-i]\cup[i,i\infty)$. We
have
\begin{equation*}
h_3(v)=4\sum_{n=0}^{\infty}\frac{1}{(2n+1)^2}\int_{0}^{v}\frac{1}{u}\bigg(\frac{u}{1+\sqrt{1+u^2}}\bigg)^{2n+1}\d
u.
\end{equation*}
The nested integrals can be evaluated by Lemma \ref{integral for
roots}. Letting $r=\frac{v}{1+\sqrt{1+v^2}}$ it is clear that
\begin{equation}\label{batir equation}
\begin{split}
h_3(v)=&4\sum_{n=0}^{\infty}\frac{1}{(2n+1)^2}\bigg(\log\left(\frac{1+r}{1-r}\right)+\frac{r^{2n+1}}{2n+1}
        -2\sum_{j=0}^{n}\frac{r^{2j+1}}{2j+1}\bigg)\\
      =&\frac{\pi^2}{2}\log\left(\frac{1+r}{1-r}\right)+4\Li_3(r)-\frac{1}{2}\Li_3(r^2)-8\F_{2,1}(1,r),
\end{split}
\end{equation}
where $\F_{2,1}(1,r)$ has a closed form provided by Eq.
\eqref{full evaluation of lambda21(x)}.  Since $|r|<1$ whenever
$v\not\in(-i\infty,-i]\cup[i,i\infty)$, we may substitute Eq.
\eqref{full evaluation of lambda21(x)} to finish the calculation.

    Observe that when $v=2$, we have
$r=\frac{\sqrt{5}-1}{2}$. It is easy to verify that
$\frac{3-\sqrt{5}}{2}=r^2=1-r=\frac{r}{1+r}$.  Using Eq. \eqref{h3
evaluated}, it follows that
\begin{equation}\label{h_3(2) formula 2}
\begin{split}
h_3(2)=&\frac{17}{2}\Li_3\left(\frac{3-\sqrt{5}}{2}\right)-4\zeta(3)
-3\log\left(\frac{1+\sqrt{5}}{2}\right)\Li_2\left(\frac{3-\sqrt{5}}{2}\right)\\
&+\frac{4\pi^2}{3}\log\left(\frac{1+\sqrt{5}}{2}\right)-\frac{26}{3}\log^3\left(\frac{1+\sqrt{5}}{2}\right).
\end{split}
\end{equation}
Eq. \eqref{h_3(2) evaluated} follows immediately from substituting
the classical formulas for
$\Li_3\left(\frac{3-\sqrt{5}}{2}\right)$ and
$\Li_2\left(\frac{3-\sqrt{5}}{2}\right)$ into Eq. \eqref{h_3(2)
formula 2}.  $\blacksquare$
\end{proof}

Notice that Eq. \eqref{h_3(2) evaluated} is equivalent to a new
evaluation of the $_4F_3$ hypergeometric function,
\begin{equation}
_4F_3\left[\substack{1,1,\frac{1}{2},\frac{1}{2}\\\frac{3}{2},\frac{3}{2},\frac{3}{2}}\bigg\vert-4\right]=\frac{7}{10}\zeta(3).
\end{equation}

\begin{corollary}\label{TS(v,1) and TS(2,1) cor}
Let $r=\frac{v}{1+\sqrt{1+v^2}}$ and suppose that
$v\not\in(-i\infty,-i]\cup[i,i\infty)$, then
\begin{align}
\begin{split}
\TS(v,1)=&\frac{\pi}{2}\int_{0}^{v}\frac{\tan^{-1}(u)}{u}\d u
        -\frac{1}{4}\Li_3(r^2)-2\Li_3(1-r)-2\Li_3\left(\frac{r}{1+r}\right)\\
        &+2\zeta(3)+\frac{1}{2}\log\left(\frac{1+r}{1-r}\right)\Li_2(r^2)+\frac{\pi^2}{3}\log(1-r)\\
        &+\frac{1}{3}\log^3(1+r)-\log(r)\log^2(1-r),
\end{split}\label{TS(v,1) closed form}\\
\TS(2,1)=&\frac{\pi}{2}\int_{0}^{2}\frac{\tan^{-1}(u)}{u}\d
u-\frac{7}{5}\zeta(3)\label{TS(2,1) evaluated}
\end{align}
\end{corollary}
\begin{proof} Eq. \eqref{TS(v,1) closed form} follows immediately
from substituting Eq. \eqref{h3 evaluated} into Eq.
\eqref{rewritten taylor series}, while  Eq. \eqref{TS(2,1)
evaluated} follows from combining Eq. \eqref{h_3(2) evaluated}
with Eq. \eqref{rewritten taylor series}.  $\blacksquare$
\end{proof}

    The fact that we can reduce $h_1(v)$, $h_2(v)$ and
$h_3(v)$ to standard polylogarithms is somewhat miraculous.
Integrating Eq. \eqref{batir equation} again, we can show that
\begin{equation}\label{h4 in terms of multi polylogs}
\begin{split}
h_4(v)=&\frac{\pi^2}{4}\left(\log(1-r^2)\log\left(\frac{1-r}{1+r}\right)+
2\Li_2\left(\frac{1-r}{2}\right)-2\Li_2\left(\frac{1+r}{2}\right)\right)\\
&+\pi^2\F_2(r)+4\F_3(r)-8\F_{3,1}(1,r)-8\F_{2,2}(1,r)\\
&+16\F_{2,1,1}(1,1,r).
\end{split}
\end{equation}
Considering the complexity of these multiple polylogarithms, it
seems unlikely that $h_n(v)$ will reduce to standard
polylogarithms for $n\ge 4$.


\section{Relations between $\S(v,1)$ and Mahler's measure, and a closed form for $\S(v,w)$.}
\label{Closed forms for S section} \init

    In this section we will study the double arcsine integral,
$\S(v,w)$.  Recall that we defined $\S(v,w)$ with an integral:
\[\S(v,w)=\int_{0}^{1}\frac{\sin^{-1}(v x)\sin^{-1}(w x)}{x}\d x.\]
First, we will show that both $\S(v,1)$ and $\S(v,v)$ reduce to
standard polylogarithms. Next, we will discuss several interesting
results relating $\S(v,1)$ and $\S(v,v)$ to Mahler's measure and
binomial sums.  Finally, Theorem \ref{S(v,w) closed form theorem}
concludes this section by expressing $\S(v,w)$ in terms of
polylogarithms.

\begin{theorem}\label{s evaluated}Assume that $0\le v\le 1$, then
$\S(v,v)$ and $\S(v,1)$ both have simple closed forms:
 \begin{align}
 \S(v,1)=&\frac{\pi}{2}\int_{0}^{v}\frac{\sin^{-1}(x)}{x}\d x-\left(\frac{\Li_3(v)-\Li_3(-v)}{2}\right),\label{closed form for
 S(v,1)}\\
 \begin{split}
\S(v,v)=&\left(\frac{\Li_3\left(e^{2i\sin^{-1}(v)}\right)+\Li_3\left(e^{-2i\sin^{-1}(v)}\right)}{4}\right)-\frac{\zeta(3)}{2}\\
&+\sin^{-1}(v)\left(\frac{\Li_2\left(e^{2i\sin^{-1}(v)}\right)-\Li_2\left(e^{-2i\sin^{-1}(v)}\right)}{2i}\right)\\
&+\left(\sin^{-1}(v)\right)^2\log(2v).
\end{split}\label{closed form for S(v,v)}
 \end{align}
\end{theorem}
\begin{proof}
To prove Eq. \eqref{closed form for S(v,1)}, we will substitute
the Taylor series for $\sin^{-1}(v x)$ into the integral
$\S(v,1)=\int_{0}^{1}\frac{\sin^{-1}(v x)\sin^{-1}(x)}{x}\d x$.
After swapping the order of summation and integration, we have
\begin{equation*}
\begin{split}
\S(v,1)&=2\sum_{n=0}^{\infty}\frac{1}{2n+1}{2n\choose
n}\left(\frac{v }{2}\right)^{2n+1}\int_{0}^{1}\sin^{-1}(x)x^{2n}\d
x\\
&=\pi\sum_{n=0}^{\infty}\frac{1}{(2n+1)^2}{2n\choose
n}\left(\frac{v}{2}\right)^{2n+1}-\sum_{n=0}^{\infty}\frac{v^{2n+1}}{(2n+1)^3}\\
&=\frac{\pi}{2}\int_{0}^{v}\frac{\sin^{-1}(x)}{x}\d
x-\left(\frac{\Li_3(v)-\Li_3(-v)}{2}\right).
\end{split}
\end{equation*}

    To prove \eqref{closed form for S(v,v)} make the $u$-substitution $x=\frac{\sin(t)}{v}$, and then
integrate by parts as follows:
\begin{equation*}
\begin{split}
\S(v,v)&=\int_{0}^{1}\frac{\left(\sin^{-1}(v x)\right)^2}{x}\d x=\int_{0}^{\sin^{-1}(v)}t^2\cot(t)\d t\\
       &=\left(\sin^{-1}(v)\right)^2\log(v)-2\int_{0}^{\sin^{-1}(v)}t\log(\sin(t))\d
       t.
\end{split}
\end{equation*}
Next substitute the Fourier series for $\log(\sin(t))$ into the
previous equation, recall that
\begin{equation*}
\log(\sin(t))=-\log(2)-\sum_{n=1}^{\infty}\frac{\cos(2n t)}{n}
\end{equation*}
is valid for $0<t< \pi$.  Integrating by parts a second time
completes the proof.  $\blacksquare$
\end{proof}

The function $\S(v,v)$ provides a connection to a second family of
interesting binomial sums.  If we recall the formula
\begin{equation*}
(\sin^{-1}(x))^2=\frac{1}{2}\sum_{n=1}^{\infty}\frac{(2x)^{2n}}{n^2{2n\choose
n}},
\end{equation*}
then it is immediately obvious that if $|v|\le 1$ we must have
\begin{equation}\label{S(v,v) taylor series}
\S(v,v)=\frac{1}{4}\sum_{n=1}^{\infty}\frac{(2v)^{2n}}{n^3{2n\choose
n}}.
\end{equation}
Comparing Eq. \eqref{S(v,v) taylor series} with Eq. \eqref{closed
form for S(v,v)} yields a classical formula:
\begin{equation}\label{S(1/2,1/2) formula sum}
\S\left(\frac{1}{2},\frac{1}{2}\right)=\frac{1}{4}\sum_{n=1}^{\infty}\frac{1}{n^3{2n\choose
n}}=\frac{1}{2}\sum_{n=1}^{\infty}\frac{\cos\left(\frac{\pi
n}{3}\right)}{n^3}-\frac{\zeta(3)}{2}+\frac{\pi}{6}\sum_{n=1}^{\infty}\frac{\sin\left(\frac{\pi
n}{3}\right)}{n^2}.
\end{equation}

\begin{proposition}\label{S related to mahler proposition}If $v\in [0,1]$ and $w\in (0,1]$, we have
\begin{equation}\label{S related to mahler}
\begin{split}
\S(v,w)=&\sin^{-1}(w)\int_{0}^{v}\frac{\sin^{-1}(u)}{u}\d
u\\
&-\frac{\pi}{2}\int_{0}^{\sin^{-1}(w)}\m\left(\frac{2v}{w}\sin(\theta)+y+z\right)\d\theta
\end{split}
\end{equation}
\end{proposition}
\begin{proof}This proof is similar to the proof of Proposition
\ref{ts related to mahler}.  After an integration by parts, and
the $u$-substitution $u=\sin(\theta)/w$, we obtain
\begin{equation*}
\S(v,w)=\sin^{-1}(w)\int_{0}^{v}\frac{\sin^{-1}(u)}{u}\d
u-\frac{\pi}{2}\int_{0}^{\sin^{-1}(w)}\int_{0}^{\frac{v}{w}\sin(\theta)}\frac{\sin^{-1}(z)}{z}\d
z\d \theta.
\end{equation*}
Since $0\le v\le 1$ and $0< w \le 1$, it follows that $0\le
\frac{v}{w}\sin(\theta)\le 1$.  Therefore we may complete the
proof by substituting Eq. \eqref{arcsine integral to mahler
measure} for the nested arcsine integral. $\blacksquare$
\end{proof}

\begin{corollary}We can recover Vandervelde's formula by
letting $w=1$ in Eq. \eqref{S related to mahler}:
\begin{equation}
\begin{split}
\m\left(v(1+x)+y+z\right)&=\frac{2}{\pi}\int_{0}^{v}\frac{\sin^{-1}(u)}{u}\d
u-\frac{4}{\pi^2}\S(v,1)\\
        &=\frac{4}{\pi^2}\left(\frac{\Li_3(v)-\Li_3(-v)}{2}\right)
\end{split}
\end{equation}
\end{corollary}

Notice that if $v=w=\frac{1}{2}$ in Eq. \eqref{S related to
mahler}, we have
\begin{align}
\S\left(\frac{1}{2},\frac{1}{2}\right)&=\frac{\pi}{6}\int_{0}^{1/2}\frac{\sin^{-1}(u)}{u}\d u-\frac{\pi}{2}\int_{0}^{\pi/6}\m\left(2\sin(\theta)+y+z\right)\d\theta\nonumber\\
                                      &=\frac{\pi}{6}\int_{0}^{1/2}\frac{\sin^{-1}(u)}{u}\d u-\frac{\pi^2}{12}\m\left(1-x^{1/6}+y+z\right)\label{S(1/2,1/2) as fractional mahler}
\end{align}
Comparing Eq. \eqref{S(1/2,1/2) as fractional mahler} to Eq.
\eqref{S(1/2,1/2) formula sum} allows us to express a famous
binomial sum as the Mahler measure of a three-variable algebraic
function.

    The final result of this section allows us to express $\S(v,w)$ in terms
of standard polylogarithms.

\begin{theorem}\label{S(v,w) closed form theorem} Suppose that $0\le v < w \le
1$, and let $\theta=\sin^{-1}(w)-\sin^{-1}(v)$.  Then we have
\begin{equation}\label{S(v,w) closed form}
\begin{split}
2\S(v,w)=&\S(v,v)+\S(w,w)-\S\left(\sin(\theta),\sin(\theta)\right)\\
&-2\Li_3\left(\frac{v}{w}\right)+\Li_3\left(\frac{v}{w}e^{i\theta}\right)+
\Li_3\left(\frac{v}{w}e^{-i\theta}\right)\\
&-i\theta\Li_2\left(\frac{v}{w}e^{i\theta}\right)+i\theta\Li_2\left(\frac{v}{w}e^{-i\theta}\right)\\
&+\frac{\theta^2}{2}\log\left(1+\frac{v^2}{w^2}-\frac{2v}{w}\cos(\theta)\right).
\end{split}
\end{equation}

    Notice that Eq. \eqref{closed form for S(v,v)} reduces $\S(v,v)$,
$\S(w,w)$, and $\S\left(\sin(\theta),\sin(\theta)\right)$ to
standard polylogarithms.
\end{theorem}
\begin{proof} The details of this proof are not particularly
difficult.  First observe the following trivial formula:
\begin{equation*}
\S(v,v)-2\S(v,w)+\S(w,w)=\int_{0}^{1}\frac{\left(\sin^{-1}(w
u)-\sin^{-1}(v u)\right)^2}{u}\d u.
\end{equation*}
Rearranging, and then applying the arcsine addition formula yields
\begin{equation}\label{s(v,w) eval intermed 1}
\begin{split}
2\S(v,w)=&\S(v,v)+\S(w,w)\\
&-\int_{0}^{1}\frac{\left(\sin^{-1}\left(w u\sqrt{1-v^2 u^2}-v
u\sqrt{1-w^2 u^2}\right)\right)^2}{u}\d u.
\end{split}
\end{equation}
This substitution is justified by the monotonicity of the arcsine
function.  In particular, $0\le v<w\le 1$ implies that
$0\le\sin^{-1}(w u)-\sin^{-1}(v u)\le \frac{\pi}{2}$ for all
$u\in[0,1]$.

    Next we will make the $u$-substitution $z=w u\sqrt{1-v^2 u^2}-v u\sqrt{1-w^2 u^2}$.  In particular,
we can show that \[u^2=\frac{z^2}{w^2+v^2-2 v w\sqrt{1-z^2}},\]
and we can easily verify that
\[\frac{1}{u}\frac{\d u}{\d z}=\frac{1}{z}-\frac{v w
z}{(v^2+w^2-2 v w\sqrt{1-z^2})\sqrt{1-z^2}}.\] Observe that the
new path of integration will run from $z=0$ to
$z=\sin(\theta)=w\sqrt{1-v^2}-v\sqrt{1-w^2}$.  Therefore, Eq.
\eqref{s(v,w) eval intermed 1} becomes
\begin{equation*}
\begin{split}
2\S(v,w)=&\S(v,v)+\S(w,w)\\
&-\int_{0}^{\sin(\theta)}\left(\sin^{-1}\left(z\right)\right)^2\left(\frac{1}{z}-\frac{v
w z}{(v^2+w^2-2 v w\sqrt{1-z^2})\sqrt{1-z^2}}\right)\d z\\
=&\S(v,v)+\S(w,w)-\S\left(\sin(\theta),\sin(\theta)\right)\\
&+\int_{0}^{\sin(\theta)}\left(\sin^{-1}\left(z\right)\right)^2\frac{v
w z}{(v^2+w^2-2 v w\sqrt{1-z^2})\sqrt{1-z^2}}\d z.
\end{split}
\end{equation*}
If we let $t=\sin^{-1}(z)$, then this last integral becomes
\begin{equation}\label{S(v,w) eval intermed 2}
\begin{split}
2\S(v,w)=&\S(v,v)+\S(w,w)-\S\left(\sin(\theta),\sin(\theta)\right)\\
                          &+\int_{0}^{\theta}t^2\frac{v w \sin(t)}{v^2+w^2-2 v w\cos(t)}\d
                          t.
\end{split}
\end{equation}

Since $0\le v <w \le 1$, a formula from \cite{Gr} shows that
\begin{equation}\label{S(v,w) eval fourier series}
\frac{vw\sin(t)}{v^2+w^2-2v
w\cos(t)}=\sum_{n=1}^{\infty}\left(\frac{v}{w}\right)^{n}\sin(n
t).
\end{equation}
The Fourier series in Eq. \eqref{S(v,w) eval fourier series}
converges uniformly since $v<w$. It follows that we may substitute
Eq. \eqref{S(v,w) eval fourier series} into Eq. \eqref{S(v,w) eval
intermed 2}, and then swap the order of summation and integration
to obtain:
\begin{equation}\label{S(v,w) eval intermed 3}
\begin{split}
2\S(v,w)=&\S(v,v)+\S(w,w)-\S\left(\sin(\theta),\sin(\theta)\right)\\
                          &+\sum_{n=1}^{\infty}\left(\frac{v}{w}\right)^n \int_{0}^{\theta}t^2\sin(n t)\d t.
\end{split}
\end{equation}
Simplifying Eq. \eqref{S(v,w) eval intermed 3} completes the proof
of Eq. \eqref{S(v,w) closed form}. $\blacksquare$
\end{proof}

\section{$q$-series for the dilogarithm, and some associated trigonometric integrals}
\label{Dilog Q series section} \init

In this section we will prove several double $q$-series expansions
for the dilogarithm. While these formulas are relatively simple,
it appears that they are new.  The first of these formulas, Eq.
\eqref{first q series for li2}, follows from a few simple
manipulations of Eq. \eqref{s evaluated}.  The remaining formulas
follow from integrals that we have evaluated in Theorem
\ref{integral table of jacobian integrals}.  Recall that Theorem
\ref{integral table of jacobian integrals} figured prominently in
the proof of Theorem \ref{arctan extra mahlers theoreom}.

    In this section, the twelve Jacobian elliptic functions will play an important
role our calculations. Recall that the Jacobian elliptic functions
are doubly periodic and meromorphic on $\mathbb{C}$. The Jacobian
sine function, $\sn(u)$, inverts the incomplete elliptical
integral of the first kind. If $u\in \mathbb{C}$ is an arbitrary
number, then under a suitable path of integration:
\begin{equation*}
u=\int_{0}^{\sn(u)}\frac{\d z}{\sqrt{(1-z^2)(1-k^2 z^2)}}.
\end{equation*}
The Jacobian amplitude can be defined by the equation
$\sn(u)=\sin(\am(u))$, and the Jacobian cosine function is defined
by $\cn(u)=\cos(\am(u))$. As usual the complementary sine function
is given by $\dn(u)=\sqrt{1-k^2\sn^2(u)}$.  Notice that every
Jacobian elliptic function implicitly depends on $k$;  this
parameter $k$ is called the elliptic modulus.

Following standard notation, we will denote the real one-quarter
period of $\sn(u)$ by $K$.  Since $\sn(K)=1$, we may compute $K$
from the usual formula
\begin{equation*}
\begin{split}
K:=K(k)&=\int_{0}^{1}\frac{\d z}{\sqrt{(1-z^2)(1-k^2z^2)}}\\
       &=\frac{\pi}{2}{_2F_1\left[\substack{\frac{1}{2},\frac{1}{2}\\1}\big|k^2\right]}.
\end{split}
\end{equation*}
Let $K'=K(\sqrt{1-k^2})$, and finally define the elliptic nome by
$q=e^{-\pi \frac{K'}{K}}$.
\begin{proposition}\label{first jacobi cn(u) integral} If $k\in (0,1)$, then we have the following
integral:
\begin{equation}\label{derivative of S(v,1) in Jacobi}
\int_{0}^{K}\am(u)\cn(u)\d
u=\frac{\pi}{2}\frac{\sin^{-1}(k)}{k}-\left(\frac{\Li_2(k)-\Li_2(-k)}{2k}\right).
\end{equation}
\end{proposition}
\begin{proof}
Taking the derivative of each side of Eq. \eqref{closed form for
S(v,1)}, we obtain:
\begin{equation}\label{derivative of S(v,1)}
\frac{\d}{\d k}\S(k,1)=\int_{0}^{1}\frac{\sin^{-1}(x)}{\sqrt{1-k^2
x^2}}\d x
=\frac{\pi}{2}\frac{\sin^{-1}(k)}{k}-\left(\frac{\Li_2(k)-\Li_2(-k)}{2k}\right).
\end{equation}
Making the $u$-substitution $x=\sn(u)$ completes the
proof.$\blacksquare$
\end{proof}

We will need the following two inversion formulas for the elliptic
nome.

\begin{lemma}\label{invert q} Let $q$ be the usual elliptic nome.  Suppose that
$q\in(0,1)$, then $q$ is invertible using either of the formulas:
\begin{align}
k&=\sin\left(4\sum_{n=0}^{\infty}\frac{(-1)^n}{2n+1}\frac{q^{n+1/2}}{(1+q^{2n+1})}\right),\label{inverse 1 of q}\\
k&=\tanh\left(4\sum_{n=0}^{\infty}\frac{1}{2n+1}\frac{q^{n+1/2}}{(1-q^{2n+1})}\right).\label{inverse
2 of q}
\end{align}
\end{lemma}
\begin{proof}
To prove Eq. \eqref{inverse 1 of q} observe that
\begin{align}
\sin^{-1}(k)&=k\int_{0}^{1}\frac{\d x}{\sqrt{1-k^2 x^2}}\notag\\
            &=k\int_{0}^{K}\cn(u)\d u\label{arcsine cn integral}
\end{align}
Recall the Fourier series expansion \cite{Gr} for $\cn(u)$:
\begin{equation}\label{cn fourier series}
\cn(u)=\frac{2\pi}{k
K}\sum_{n=0}^{\infty}\frac{q^{n+1/2}}{1+q^{2n+1}}\cos\left(\frac{\pi(2n+1)}{2K}u\right).
\end{equation}
Since $0<q<1$, this Fourier series converges uniformly. It follows
that we may substitute Eq. \eqref{cn fourier series} into Eq.
\eqref{arcsine cn integral}, and then swap the order of summation
and integration to obtain:
\begin{align}
\sin^{-1}(k)&=\frac{2\pi}{
K}\sum_{n=0}^{\infty}\frac{q^{n+1/2}}{1+q^{2n+1}}\int_{0}^{K}\cos\left(\frac{\pi(2n+1)}{2K}u\right)\d
u\notag\\
    &=4\sum_{n=0}^{\infty}\frac{(-1)^n}{2n+1}\frac{q^{n+1/2}}{(1+q^{2n+1})}.
\end{align}
Eq. \eqref{inverse 1 of q} follows immediately from taking the
sine of both sides of the equation.

Eq. \eqref{inverse 2 of q} can be proved in a similar manner when
starting from the integral \[\tanh^{-1}(k)=k\int_{0}^{1}\frac{\d
x}{1-k^2 x^2}.\] $\blacksquare$
\end{proof}

Next we will utilize the Fourier-series expansions for the
Jacobian elliptic functions to prove the following theorem:

\begin{theorem}\label{first q series for li2 theorem} If $q$ is the usual elliptic nome, then the following formula holds for the dilogarithm:
\begin{equation}\label{first q series for li2}
\begin{split}
\frac{\Li_2(k)-\Li_2(-k)}{8}=&\sum_{n=0}^{\infty}\frac{1}{(2n+1)^2}
\frac{q^{n+1/2}}{(1+q^{2n+1})}\\
&+4\sum_{\substack{n=0\\m=1}}^{\infty}\frac{1}{(2n+1)^2-(2m)^2}\frac{q^{n+m+1/2}}{(1+q^{2m})(1+q^{2n+1})}
\end{split}
\end{equation}
\end{theorem}
\begin{proof}
We have already stated the Fourier series expansion for $\cn(u)$
in Eq. \eqref{cn fourier series}.  We will also require the
Fourier series \cite{Gr} for $\am(u)$:
\begin{equation}\label{am fourier series}
\am(u)=\frac{\pi}{2K}u+2\sum_{n=1}^{\infty}\frac{1}{n}\frac{q^n}{1+q^{2n}}\sin\left(\frac{\pi
n}{K}u\right).
\end{equation}
Substituting Eq. \eqref{cn fourier series} and Eq. \eqref{am
fourier series} into the integral in Eq. \eqref{derivative of
S(v,1) in Jacobi}, and then simplifying yields:
\begin{equation}\label{polylog pre formula}
\begin{split}
\frac{\Li_2(k)-\Li_2(-k)}{8}&-\frac{\pi}{8}\sin^{-1}(k)\\
=-&\frac{\pi}{2}\sum_{n=0}^{\infty}\frac{(-1)^n}{2n+1}
\frac{q^{n+1/2}}{(1+q^{2n+1})}+\sum_{n=0}^{\infty}\frac{1}{(2n+1)^2}
\frac{q^{n+1/2}}{(1+q^{2n+1})}\\
&+4\sum_{\substack{n=0\\m=1}}^{\infty}\frac{1}{(2n+1)^2-(2m)^2}\frac{q^{n+m+1/2}}{(1+q^{2m})(1+q^{2n+1})}.
\end{split}
\end{equation}

    This proof is nearly complete, the final step is to substitute the
identity
\begin{equation*}
\sin^{-1}(k)=4\sum_{n=0}^{\infty}\frac{(-1)^n}{2n+1}\frac{q^{n+1/2}}{(1+q^{2n+1})}
\end{equation*}
into Eq. \eqref{polylog pre formula}.  This formula for
$\sin^{-1}(k)$ follows immediately from Lemma \ref{invert q}.
$\blacksquare$
\end{proof}

The fact that Eq. \eqref{first q series for li2} follow easily
from an integral of the form
\begin{equation*}
\int_{0}^{K}\am(u)\varphi(u)\d u,
\end{equation*}
suggests that we should try to generalize Eq. \eqref{first q
series for li2} by allowing $\varphi(u)$ to equal one of the other
eleven Jacobian elliptic functions.  Theorem \ref{integral table
of jacobian integrals} proves that ten of these eleven integrals
reduce to dilogarithms and elementary functions. First, Theorem
\ref{evaluate mahlers measure} will prove that the one exceptional
integral can be expressed as the Mahler measure of an elliptic
curve.

\begin{theorem}\label{evaluate mahlers measure}
The following formulas hold whenever $k\in (0,1]$:
\begin{align}
\m\left(\frac{4}{k}+x+\frac{1}{x}\right.&\left.+y+\frac{1}{y}\right)\nonumber\\
&=-\log\left(\frac{k}{1+\sqrt{1-k^2}}\right)
+\frac{2}{\pi}\int_{0}^{1}\frac{\sin^{-1}(x)}{x\sqrt{1-k^2 x^2}}\d
x\label{elliptic mahler formula 1}\\
&=-\log\left(\frac{k}{1+\sqrt{1-k^2}}\right)+\frac{2}{\pi}\int_{0}^{K}\am(u)\frac{\cn(u)}{\sn(u)}\d
u\label{elliptic mahler formula 2}.
\end{align}
\end{theorem}
\begin{proof}
First observe that if $k\in\mathbb{R}$ and $0<k\le1$, then
\begin{equation*}
\m\left(\frac{4}{k}+x+\frac{1}{x}+y+\frac{1}{y}\right)=-\log\left(\frac{k}{4}\right)+\m\left(1+\frac{k}{4}\left(x+\frac{1}
{x}+y+\frac{1}{y}\right)\right).
\end{equation*}
For brevity let $
\varphi(k)=\m\left(1+\frac{k}{4}\left(x+\frac{1}{x}+y+\frac{1}{y}\right)\right)$.
Making the change of variables $(x,y)\rightarrow(x/y,yx)$, we
have:
\begin{align*}
\varphi(k)=&\m\left(1+\frac{k}{4}\left(x+x^{-1}\right)\left(y+y^{-1}\right)\right)\\
          =&\m\left(\frac{k}{4}\left(y+y^{-1}\right)\right)+\m\left(x^2+\frac{4}{k}\left(\frac{1}{y+y^{-1}}\right)x+1\right)\\
          =&\log\left(\frac{k}{4}\right)+\m\left(x^2+\frac{4}{k}\left(\frac{1}{y+y^{-1}}\right)x+1\right).
\end{align*}
Applying Jensen's formula with respect to $x$ reduces $\varphi(k)$
to a pair of one-dimensional integrals:
\begin{equation}\label{elliptic mahler integral log +}
\begin{split}
\varphi(k)=\log\left(\frac{k}{4}\right)+&\frac{1}{2\pi}\int_{0}^{2\pi}\log^{+}\left\vert\frac{1+\sqrt{1-k^2\cos^2(\theta)}}{k\cos(\theta)}\right\vert\d\theta\\
&+\frac{1}{2\pi}\int_{0}^{2\pi}\log^{+}\left\vert\frac{1-\sqrt{1-k^2\cos^2(\theta)}}{k\cos(\theta)}\right\vert\d\theta.
\end{split}
\end{equation}
The right-hand integral vanishes under the assumption that $0<k\le
1$. Therefore, it follows that Eq. \eqref{elliptic mahler integral
log +} reduces to
\begin{equation*}
\varphi(k)=\log\left(\frac{k}{4}\right)+
\frac{2}{\pi}\int_{0}^{\pi/2}\log\left(\frac{1+\sqrt{1-k^2\cos^2(\theta)}}{k\cos(\theta)}\right)\d\theta.
\end{equation*}
With the observation that
$\int_{0}^{\pi/2}\log\left(\cos(\theta)\right)\d
\theta=-\frac{\pi}{2}\log(2)$, this formula becomes:
\begin{equation}
\varphi(k)=\frac{2}{\pi}\int_{0}^{\pi/2}\log\left(\frac{1+\sqrt{1-k^2\cos^2(\theta)}}{2}\right)\d\theta.
\end{equation}
Making the $u$-substitution of $x=\cos(\theta)$, we obtain
\begin{equation*}
\varphi(k)=\frac{2}{\pi}\int_{0}^{1}\log\left(\frac{1+\sqrt{1-k^2
x^2}}{2}\right)\frac{1}{\sqrt{1-x^2}}\d x.
\end{equation*}
Integrating by parts to eliminate the logarithmic term yields:
\begin{equation*}
\begin{split}
\varphi(k)&=\log\left(\frac{1+\sqrt{1-k^2}}{2}\right)+\frac{2}{\pi}
            \int_{0}^{1}\frac{\sin^{-1}(x)}{x}\left(\frac{1-\sqrt{1-k^2
            x^2}}{\sqrt{1-k^2 x^2}}\right)\d x\\
        &=\log\left(\frac{1+\sqrt{1-k^2}}{2}\right)+\frac{2}{\pi}
            \int_{0}^{1}\frac{\sin^{-1}(x)}{x\sqrt{1-k^2 x^2}}\d
            x-\frac{2}{\pi} \int_{0}^{1}\frac{\sin^{-1}(x)}{x}\d x.
\end{split}
\end{equation*}
Since $\int_{0}^{1}\frac{\sin^{-1}(x)}{x}\d
x=\frac{\pi}{2}\log(2)$, it follows that
\begin{equation*}
\varphi(k)=\log\left(\frac{1+\sqrt{1-k^2}}{4}\right)+\frac{2}{\pi}\int_{0}^{1}\frac{\sin^{-1}(x)}{x\sqrt{1-k^2
x^2}}\d x,\end{equation*}
 from which we obtain
\begin{equation*}
\m\left(\frac{4}{k}+x+\frac{1}{x}+y+\frac{1}{y}\right)=-\log\left(\frac{k}{1+\sqrt{1-k^2}}\right)
+\frac{2}{\pi}\int_{0}^{1}\frac{\sin^{-1}(x)}{x\sqrt{1-k^2 x^2}}\d x.\\
\end{equation*}

To prove Eq. \eqref{elliptic mahler formula 2} simply make the
$u$-substitution $x=\sn(u)$.$\blacksquare$
\end{proof}

The elliptic curve defined by the equation $4/k+x+1/x+y+1/y=0$ was
one of the simplest curves that Boyd studied in \cite{Bo1}.
Rodriguez Villegas derived $q$-series expansions for a wide class
of functions defined by the Mahler measures of elliptic curves in
\cite{RV}. We can recover one of his results by substituting the
Fourier series expansions for $\am(u)$ and $\cn(u)/\sn(u)$ into
Eq. \eqref{elliptic mahler formula 2}.

    If we let $k=\sin(\theta)$, and then integrate
Eq. \eqref{elliptic mahler formula 1} from $\theta=0$ to
$\theta=\frac{\pi}{2}$, we can prove that
\begin{equation}
\m\left(8+\left(z+\frac{1}{z}\right)\left(x+
\frac{1}{x}+y+\frac{1}{y}\right)\right)=\frac{4}{\pi}G+\frac{4}{\pi^2}\int_{0}^{1}\frac{\sin^{-1}(x)}{x}K(x)\d
x.
\end{equation}
Using Mathematica, we can reduce the right-hand integral to a
rather complicated expression involving balanced hypergeometric
functions evaluated at one.

\begin{theorem}\label{integral table of jacobian integrals} We will assume that $0< k< 1$ and that each
Jacobian elliptic function has modulus $k$.  Let
$p=\sqrt{\frac{1-k}{1+k}}$, $r=\frac{k}{1+\sqrt{1-k^2}}$, and
$s=\frac{k}{\sqrt{1-k^2}}$, then {\allowdisplaybreaks
\begin{align}
\begin{split}
\int_{0}^{K}\am(u)\sn(u)\d u
&=\int_{0}^{1}\frac{u\sin^{-1}(u)}{\sqrt{(1-u^2)(1-k^2 u^2)}}\d
u\\
 &=\frac{\Li_2(i s)-\Li_2(-i s)}{2k i}\label{sn(u) integral}\end{split}\\
\begin{split}\int_{0}^{K}\am(u)\cn(u)\d
u&=\int_{0}^{1}\frac{\sin^{-1}(u)}{\sqrt{1-k^2 u^2}}\d u\\
 &=\frac{\pi}{2}\frac{\sin^{-1}(k)}{k}-\frac{\Li_2(k)-\Li_2(-k)}{2k}\label{cn(u) integral}\end{split}\\
\int_{0}^{K}\am(u)\dn(u)\d u &=\frac{\pi^2}{8}\label{dn(u) integral}\\
\begin{split}\int_{0}^{K}\am(u)\frac{1}{\sn(u)}\d
u&=\int_{0}^{1}\frac{\sin^{-1}(u)}{u\sqrt{(1-u^2)(1-k^2 u^2)}}\d
u\\
&=-\frac{\pi}{2}\log\left(p\right)+\frac{\Li_2(i
 p)-\Li_2(-i p)}{i}\label{1/sn(u) integral}\end{split}\\
\int_{0}^{K}\am(u)\frac{1}{\cn(u)}\d u&=\infty\label{1/cn(u) integral}\\
\begin{split}\int_{0}^{K}\am(u)\frac{1}{\dn(u)}\d
u&=\int_{0}^{1}\frac{\sin^{-1}(u)}{(1-k^2 u^2)\sqrt{1-u^2}}\d u\\
  &=\frac{1}{\sqrt{1-k^2}}\left(\frac{\pi^2}{8}+\frac{\Li_2(r^2)-\Li_2(-r^2)}{2}\right)\label{1/dn(u) integral}\end{split}\\
\int_{0}^{K}\am(u)\frac{\sn(u)}{\cn(u)}\d
u&=\infty\label{sn(u)/cn(u)
integral}\\
\begin{split}\int_{0}^{K}\am(u)\frac{\sn(u)}{\dn(u)}\d
u&=\int_{0}^{1}\frac{u\sin^{-1}(u)}{(1-k^2 u^2)\sqrt{1-u^2}}\d u\\
 &=\frac{\Li_2(r)-\Li_2(-r)}{k\sqrt{1-k^2}}\end{split}\label{sn(u)/dn(u) integral}\\
\begin{split}\int_{0}^{K}\am(u)\frac{\cn(u)}{\sn(u)}\d
u&=\int_{0}^{1}\frac{\sin^{-1}(u)}{u\sqrt{1-k^2 u^2}}\d u\\
 &=\frac{\pi}{2}\log\left(r\right)+\frac{\pi}{2}\m\left(\frac{4}{k}+x+\frac{1}{x}+y+\frac{1}{y}\right)\label{cn(u)/sn(u) integral}\end{split}\\
\begin{split}\int_{0}^{K}\am(u)\frac{\cn(u)}{\dn(u)}\d u &=
\int_{0}^{1}\frac{\sin^{-1}(u)}{1-k^2 u^2}\d u\\
&=-\frac{\pi}{2k}\log(p)-\frac{\Li_2(i
r)-\Li_2(-i r)}{k i}\label{cn(u)/dn(u) integral}\end{split}\\
\int_{0}^{K}\am(u)\frac{\dn(u)}{\sn(u)}\d u &=2G\label{dn(u)/sn(u) integral}\\
\int_{0}^{K}\am(u)\frac{\dn(u)}{\cn(u)}\d u
&=\infty\label{dn(u)/cn(u) integral}
\end{align}}
\end{theorem}
\begin{proof}  First observe that Eq. \eqref{1/cn(u) integral}, Eq. \eqref{sn(u)/cn(u) integral}, and
Eq. \eqref{dn(u)/cn(u) integral} all follow from the fact that
$\cn(K)=0$.  Similarly, Eq. \eqref{dn(u) integral} and Eq.
\eqref{dn(u)/sn(u) integral} both follow from the formula
$\frac{\d}{\d u}\am(u)=\dn(u)$.

    We already proved Eq. \eqref{cn(u) integral} in Proposition \ref{first jacobi cn(u) integral},
 and Eq. \eqref{cn(u)/sn(u) integral} was proved in Theorem \ref{evaluate mahlers measure}. This leaves
a total of five formulas to prove.

To prove Eq. \eqref{sn(u) integral}, observe that after letting
$u=\sqrt{1-z^2}$, we have
\begin{equation*}
\int_{0}^{1}\frac{u\sin^{-1}(u)}{\sqrt{(1-u^2)(1-k^2 u^2)}}\d
u=\frac{i s}{k
i}\int_{0}^{1}\frac{\sin^{-1}(\sqrt{1-z^2})}{\sqrt{1-(i s)^2
z^2}}\d z.
\end{equation*}
If $0< k\le 1/\sqrt{2}$, then $|s|\le1$. With this restriction on
$k$, we may expand the square root in a Taylor series to obtain:
\begin{align}
&=\frac{1}{k i}\sum_{m=0}^{\infty}(-1)^m{-1/2\choose
m}(is)^{2m+1}\int_{0}^{1}\sin^{-1}(\sqrt{1-z^2})z^{2m}\d z\notag\\
&=\frac{1}{k i}\sum_{m=0}^{\infty}\frac{(i
s)^{2m+1}}{(2m+1)^2}\notag\\
&=\frac{\Li_2(is)-\Li_2(-i s)}{2ki}\label{dilog temp sum}.
\end{align}
Notice that Eq. \eqref{dilog temp sum} extends to $0< k < 1$,
since both sides of the equation are analytic in this interval.
Therefore, Eq. \eqref{sn(u) integral} follows immediately.

    To prove Eq. \eqref{1/sn(u) integral} make the
$u$-substitution $u=\frac{z}{\sqrt{1-k^2+z^2}}$. Recalling that $
\sin^{-1}\left(\frac{z}{\sqrt{1-k^2+z^2}}\right)=\tan^{-1}\left(\frac{z}{\sqrt{1-k^2}}\right)$,
we obtain
\begin{equation*}
\int_{0}^{1}\frac{\sin^{-1}(u)}{u\sqrt{(1-u^2)(1-k^2 u^2)}}\d
u=\int_{0}^{\infty}\frac{\tan^{-1}\left(\frac{z}{\sqrt{1-k^2}}\right)}{z\sqrt{1+z^2}}\d
z.
\end{equation*}
Using Mathematica to evaluate this last integral yields:
\begin{align*}
&=-\frac{\pi}{2}\log(p)+\sqrt{1-k^2}{_3F_2\left[\substack{\frac{1}{2},1,1\\\frac{3}{2},\frac{3}{2}}\bigg{\vert}
1-k^2\right]}\\
&=-\frac{\pi}{2}\log(p)+\frac{1}{2}\sum_{n=0}^{\infty}\frac{\left(2\sqrt{1-k^2}\right)^{2n+1}}{(2n+1)^2{2n\choose
n}}\\
&=-\frac{\pi}{2}\log(p)+2\left(\frac{\Li_2(ip)-\Li_2(-ip)}{2i}\right),
\end{align*}
where Eq. \eqref{h 2 identity} justifies the final step.

    To prove Eq. \eqref{1/dn(u) integral} observe that after the
$u$-substitution $u=\sin(\theta)$ we have
\begin{equation*}
\int_{0}^{1}\frac{\sin^{-1}(u)}{(1-k^2u^2)\sqrt{1-u^2}}\d
u=\int_{0}^{\pi/2}\frac{\theta}{1-k^2\sin^2(\theta)}\d\theta.
\end{equation*}
Now substitute the Fourier series
\begin{equation}\label{fourier series 2}
\frac{\sqrt{1-k^2}}{1-k^2\sin^2(\theta)}=1+2\sum_{m=1}^{\infty}(-1)^m\left(\frac{k}{1+\sqrt{1-k^2}}
\right)^{2m}\cos(2m\theta)
\end{equation}
into the integral, and simplify to complete the proof.

    The proof of Eq. \eqref{sn(u)/dn(u) integral} follows
the same lines as the derivation of Eq. \eqref{1/dn(u) integral}.
Observe that
\begin{equation*}
\int_{0}^{1}\frac{u\sin^{-1}(u)}{(1-k^2 u^2)\sqrt{1-u^2}}\d
u=\int_{0}^{\pi/2}\frac{\theta\sin(\theta)}{1-k^2\sin^2(\theta)}\d\theta.
\end{equation*}
Now substitute the Fourier series
\begin{equation}
\frac{k\sqrt{1-k^2}\sin(\theta)}{1-k^2\sin^2(\theta)}=2\sum_{m=0}^{\infty}(-1)^m\left(\frac{k}{1+\sqrt{1-k^2}}\right)^{2m+1}
\sin\left((2m+1)\theta\right)
\end{equation}
into the integral, and simplify to complete the proof.

    Finally, we are left with Eq. \eqref{cn(u)/dn(u) integral}.
Expanding $1/(1-k^2 u^2)$ in a geometric series yields:
\begin{align*}
\int_{0}^{1}\frac{\sin^{-1}(u)}{1-k^2 u^2}\d u&=
\sum_{n=0}^{\infty}k^{2n}\int_{0}^{1}\sin^{-1}(u)u^{2n}\d u\\
&=\sum_{n=0}^{\infty}k^{2n}\left(\frac{\pi/2}{2n+1}-\frac{2^{2n}}{(2n+1)^2{2n\choose
n}}\right)\\
&=-\frac{\pi}{4k}\log\left(\frac{1-k}{1+k}\right)-\frac{h_2(i
k)}{2i k}.
\end{align*}
Substituting the closed form for $h_2(ik)$ provided by Eq.
\eqref{h 2 identity} completes the proof. $\blacksquare$
\end{proof}

We can obtain each of the following $q$-series by applying the
method from Theorem \ref{first q series for li2 theorem} to the
formulas in Theorem \ref{integral table of jacobian integrals}.
\begin{corollary}\label{five q series for dilog} Let $p=\sqrt{\frac{1-k}{1+k}}$, and let $r=\frac{k}{1+\sqrt{1-k^2}}$.
The following formulas hold for the
dilogarithm:{\allowdisplaybreaks
\begin{align}
\begin{split}
\frac{\Li_2(k)-\Li_2(-k)}{8}=&\sum_{n=0}^{\infty}\frac{1}{(2n+1)^2}
\frac{q^{n+1/2}}{(1+q^{2n+1})}\\
&+4\sum_{\substack{n=0\\m=1}}^{\infty}\frac{1}{(2n+1)^2-(2m)^2}\frac{q^{m+n+1/2}}{(1+q^{2m})(1+q^{2n+1})},
\end{split}\label{first dilog q-series}\\
\begin{split}
\frac{\Li_2(r)-\Li_2(-r)}{4}=&\sum_{n=0}^{\infty}\frac{1}{(2n+1)^2}
\frac{q^{n+1/2}}{(1+q^{2n+1})}\\
&+4\sum_{\substack{n=0\\m=1}}^{\infty}\frac{(-1)^m}{(2n+1)^2-(2m)^2}\frac{q^{m+n+1/2}}{(1+q^{2m})(1+q^{2n+1})},
\end{split}\label{second dilog q-series}\\
\begin{split}
\frac{\Li_2(i p)-\Li_2(-i
p)}{8i}=&\frac{G}{4}+\frac{\pi}{16}\log(p)+\sum_{n=0}^{\infty}\frac{(-1)^n}{(2n+1)^2}
\frac{q^{2n+1}}{(1-q^{4n+2})}\\
&+4\sum_{\substack{n=0\\m=1}}^{\infty}\frac{(-1)^{n+m}}{(2n+1)^2-(2m)^2}\frac{q^{m+2n+1}}{(1+q^{2m})(1-q^{4n+2})}.
\end{split}\label{third dilog q series}
\end{align}}
\end{corollary}
\begin{proof}  As we have already stated, each of these formulas
can be proved by substituting Fourier series expansions for the
Jacobian elliptic functions into Theorem \ref{integral table of
jacobian integrals}.

    Using the method described, we have already
proved Eq. \eqref{first dilog q-series} in Theorem \ref{first q
series for li2 theorem}. Eq. \eqref{second dilog q-series} follows
in a similar manner from Eq. \eqref{sn(u)/dn(u) integral}.

    Eq. \eqref{third dilog q series} is a little trickier to prove.
Expanding Eq. \eqref{dn(u)/sn(u) integral} in a $q$-series yields
the identity
\begin{equation}\label{auxiliary q-series}
\begin{split}
\sum_{n=1}^{\infty}\frac{1}{n}\frac{q^n}{1+q^{2n}}&\sum_{j=0}^{n-1}\frac{(-1)^j}{2j+1}\\
=&\sum_{n=0}^{\infty}\frac{(-1)^n}{(2n+1)^2}\frac{q^{2n+1}}{1+q^{2n+1}}\\
&+4\sum_{\substack{n=0\\m=1}}^{\infty}\frac{(-1)^{n+m}}{(2n+1)^2-(2m)^2}\frac{q^{m+2n+1}}{(1+q^{2m})(1+q^{2n+1})}.
\end{split}
\end{equation}
Next expand Eq. \eqref{1/sn(u) integral} in the $q$-series
\begin{equation*}
\begin{split}
\frac{\Li_2(i p)-\Li_2(-i
p)}{4i}=&\frac{G}{2}+\frac{\pi}{8}\log(p)\\
        &+\sum_{n=1}^{\infty}\frac{1}{n}\frac{q^n}{1+q^{2n}}\sum_{j=0}^{n-1}\frac{(-1)^j}{2j+1}
+\sum_{n=0}^{\infty}\frac{(-1)^n}{(2n+1)^2}\frac{q^{2n+1}}{1-q^{2n+1}}\\
&+4\sum_{\substack{n=0\\m=1}}^{\infty}\frac{(-1)^{n+m}}{(2n+1)^2-(2m)^2}\frac{q^{m+2n+1}}{(1+q^{2m})(1-q^{2n+1})},
\end{split}
\end{equation*}
and then combine it with Eq. \eqref{auxiliary q-series} to
complete the proof of Eq. \eqref{third dilog q
series}.$\blacksquare$
\end{proof}

It is important to notice that the nine convergent integrals in
Theorem \ref{integral table of jacobian integrals} only produce
three interesting $q$-series for the dilogarithm.  The other
$q$-series we may obtain from Theorem \ref{integral table of
jacobian integrals} really just restate known facts about the
elliptic nome. For example, if we expand Eq. \eqref{1/dn(u)
integral} in a $q$-series, we will obtain Eq. \eqref{first dilog
q-series} with $q$ replaced by $q^2$ and $k$ replaced by $r^2$.
This is equivalent to the fact that
$q\left(\left(\frac{k}{1+\sqrt{1-k^2}}\right)^2\right)=q^2(k)$. If
we let $\ell=\left(\frac{k}{1+\sqrt{1-k^2}}\right)^2$, then
clearly $k$ and $\ell$ satisfy a second degree modular equation
\cite{Be}.
\section{A closed form for $\T(v,w)$, and Mahler measures for
$\T\left(v,\frac{1}{v}\right)$} \label{Closed form for T section}
\init
Recall that we defined $\T(v,w)$ using the following integral:
\begin{equation}\label{T(v,w) definition}
\T(v,w)=\int_{0}^{1}\frac{\tan^{-1}(v x)\tan^{-1}(w x)}{x}\d x.
\end{equation}
Since this integral involves two arctangents, rather than one or
two arcsines, $\T(v,w)$ possesses a number of useful properties
that $\S(v,w)$ and $\TS(v,w)$ appear to lack.

    First observe that $\T(v,w)$ obeys an eight term
functional equation.  If we let
$\T(v)=\int_{0}^{v}\frac{\tan^{-1}(x)}{x}\d x$, then we can use
properties of the arctangent function to prove the following
formula:
\begin{equation}\label{T(v,w) functional equation}
\begin{split}
\T(v,w)+\T\left(\frac{1}{v},\frac{1}{w}\right)&-\T\left(\frac{w}{v},1\right)-\T\left(\frac{v}{w},1\right)\\
&=\frac{\pi}{2}\left(\T(v)+\T\left(\frac{1}{w}\right)-\T\left(\frac{v}{w}\right)-\T\left(1\right)\right).
\end{split}
\end{equation}

    If $|v|<1$ and $|w|<1$, we can substitute arctangent Taylor series
expansions into Eq. \eqref{T(v,w) definition} to obtain:
\begin{equation}\label{T(v,w) series expansion}
\begin{split}
\T(v,w)=&\sum_{n=0}^{\infty}\frac{(-1)^n
w^{2n+2}}{(2n+2)^2}\sum_{m=0}^{n}\frac{\left(v/w\right)^{2m+1}}{2m+1}\\
&+\sum_{n=0}^{\infty}\frac{(-1)^n
v^{2n+2}}{(2n+2)^2}\sum_{m=0}^{n}\frac{\left(w/v\right)^{2m+1}}{2m+1}.
\end{split}
\end{equation}
Eq. \eqref{T(v,w) series expansion} immediately reduces $\T(v,w)$
to multiple polylogarithms.  Theorem \ref{T(v,w) closed form
theorem} improves upon this result by expressing $\T(v,w)$ in
terms of standard polylogarithms.
\begin{theorem}\label{T(v,w) closed form theorem} If $v$ and $w$ are real numbers such that $|w/v|\le 1$, then
{\allowdisplaybreaks
\begin{equation}\label{T(v,w) closed form}
\begin{split}
 -4 \T(v,w)=&
 2\Li_3\left(\frac{w}{v}\right)-2\Li_3\left(-\frac{w}{v}\right)
          +\Li_3\left(\frac{1-v i}{1-wi}\right)+\Li_3\left(\frac{1+v i}{1+w
          i}\right)\\
          &-\Li_3\left(\frac{1+v i}{1-w i}\right)-\Li_3\left(\frac{1-v i}{1+w
            i}\right)\\
          &-\Li_3\left(\frac{w(1-v i)}{v(1-w i)}\right)-\Li_3\left(\frac{w(1+v i)}{v(1+w
            i)}\right)\\
          &+\Li_3\left(-\frac{w(1+v i)}{v(1-w
            i)}\right)+\Li_3\left(-\frac{w(1-vi)}{v(1+w i)}\right)\\
          &+\log\left(\frac{1+v^2}{1+w^2}\right)\left(\Li_2\left(\frac{w}{v}\right)-\Li_2\left(-\frac{w}{v}\right)\right)\\
          &-4\tan^{-1}(v)\left(\frac{\Li_2(w i)-\Li_2(-w i)}{2i}\right)\\
          &-4\tan^{-1}(w)\left(\frac{\Li_2(v i)-\Li_2(-v i)}{2 i}\right)\\
          &-\pi
          \log\left(\frac{1+v^2}{1+w^2}\right)\tan^{-1}(w)+4\log(v)\tan^{-1}(v)\tan^{-1}(w).
\end{split}
\end{equation}}
\end{theorem}
\begin{proof}
Substituting logarithms for the inverse tangents, we obtain
\begin{equation*}
\begin{split}
-4\T(v,w)=&\int_{0}^{1}\log\left(\frac{1+i v u}{1-i
v u}\right)\log\left(\frac{1+i w u}{1-i w u}\right)\frac{\d u}{u}\\
      =&\int_{0}^{i
      w}\log\left(\frac{1+\frac{v}{w}u}{1-\frac{v}{w}u}\right)\log\left(\frac{1+u}{1-u}\right)\frac{\d
      u}{u}.
\end{split}
\end{equation*}
The identity then follows (more or less) immediately from four
applications of Lewin's formula
\begin{equation}\label{Lewin's Formula}
\begin{split}
\int_{0}^{x}\log\left(1-z\right)&\log\left(1-c z\right)\frac{\d
z}{z}\\
=&\Li_3\left(\frac{1-c
x}{1-x}\right)+\Li_3\left(\frac{1}{c}\right)+\Li_3(1)\\
&-\Li_3(1-c x)-\Li_3(1-x)-\Li_3\left(\frac{1-c
x}{c(1-x)}\right)\\
&+\log(1-c
x)\left[\Li_2\left(\frac{1}{c}\right)-\Li_2(x)\right]\\
&+\log(1-x)\left[\Li_2(1-c
x)-\Li_2\left(\frac{1}{c}\right)+\frac{\pi^2}{6}\right]\\
&+\frac{1}{2}\log(c)\log^2(1-x),
\end{split}
\end{equation}
which was proved in \cite{Le1}.  Condon has discussed the
intricacies of applying this equation in \cite{Co}.$\blacksquare$
\end{proof}

    This closed form for $\T(v,w)$ is quite complicated.  Notice
that a slight change in the integrand in Eq. \eqref{T(v,w)
definition} produces a remarkably simplified formula:
\begin{equation}\label{double arctan simple integral}
\int_{0}^{1}\frac{\tan^{-1}(v x)\tan^{-1}(w x)}{\sqrt{1-x^2}}\d
x=\pi\sum_{n=0}^{\infty}\frac{\left(\frac{v}{1+\sqrt{1+v^2}}\frac{w}{1+\sqrt{1+w^2}}\right)^{2n+1}}{(2n+1)^2}.
\end{equation}
To prove Eq. \eqref{double arctan simple integral}, make the
$u$-substitution $x=\sin(\theta)$, and then apply Eq.
\eqref{arctan fourier series} twice.

    There are two special cases of Eq. \eqref{T(v,w) closed form} worth mentioning.
First observe that $\T\left(v,\frac{1}{v}\right)$ reduces to a
very simple expression. If we let $w\rightarrow 1/v$ in Eq.
\eqref{T(v,w) closed form}, and perform a few torturous
manipulations, we can show that
\begin{equation}\label{T(v,1/v) closed form}
\begin{split}
\T\left(v,\frac{1}{v}\right)=&\frac{\pi}{2}\Im\left[\Li_2(i
v)\right]-\frac{1}{2}\left(\Li_3(v^2)-\Li_3(-v^2)\right)\\
&+\frac{\log(v)}{2}\left(\Li_2(v^2)-\Li_2(-v^2)\right).
\end{split}
\end{equation}
Lal\'in obtained an equivalent form of Eq. \eqref{T(v,1/v) closed
form} using a different method. (See Appendix 2 in \cite{La1}.
Lal\'in's formula for $\T(v,1)+\T(1/v,1)$ reduces to Eq.
\eqref{T(v,1/v) closed form} after applying Eq. \eqref{T(v,w)
functional equation} with $w=1/v$).  Observe that when $w=v$ in
Eq. \eqref{T(v,w) closed form}, we have
\begin{equation}
\begin{split}
\T(v,v)=&\frac{1}{2}\Re\left[\Li_3\left(\frac{1+v i}{1-v
i}\right)-\Li_3\left(-\frac{1+v i}{1-v
i}\right)\right]-\frac{7}{8}\zeta(3)\\
&+2\tan^{-1}(v)\Im\left[\Li_2(i
v)\right]-\log(v)\left(\tan^{-1}(v)\right)^2.
\end{split}
\end{equation}
Finally, it appears that $\T(v,1)$ does not reduce to any
particularly simple expression.  Letting $w\rightarrow 1$ fails to
simplify Eq. \eqref{T(v,w) closed form} in any appreciable way.
Expanding $\T(v,1)$ in a Taylor series results in an equally
complicated expression:
\begin{equation}
\begin{split}
\T(v,1)=&\frac{1}{2}\sum_{n=0}^{\infty}\frac{v^{2n+1}}{(2n+1)^2}\sum_{k=1}^{n}\frac{(-1)^{k+1}}{k}\\
&+\frac{\pi}{4}\int_{0}^{v}\frac{\tan^{-1}(x)}{x}\d
x-\frac{\log(2)}{4}\left(\Li_2(v)-\Li_2(-v)\right)
\end{split}
\end{equation}

     Theorem \ref{T(v,1/v) mahler theorem} relates $\T(v,w)$ to three-variable Mahler measures, and
generalizes one of Lal\'in's formulas.  Once again, we will need a
simple lemma before we prove our theorem.

\begin{lemma} Suppose that $v$ and $w$ are positive real numbers,
then
\begin{equation}\label{T(v,w) as a double integral}
\T(v,w)=\tan^{-1}(v)\int_{0}^{w}\frac{\tan^{-1}(u)}{u}\d
u-\int_{0}^{\tan^{-1}(v)}\int_{0}^{\frac{w}{v}\tan(\theta)}\frac{\tan^{-1}(z)}{z}\d
z\d\theta,
\end{equation}
\begin{equation}\label{T(v,1/v) as a double integral}
\T\left(v,\frac{1}{v}\right)=\frac{\pi}{2}\int_{0}^{v}\frac{\tan^{-1}(u)}{u}\d
u-\frac{1}{2}\int_{0}^{\pi/2}\int_{0}^{v^2\tan(\theta)}\frac{\tan^{-1}(z)}{z}\d
z\d\theta.
\end{equation}
\end{lemma}
\begin{proof} While we can verify Eq. \eqref{T(v,w) as a double integral} with a trivial integration
by parts, the proof of Eq. \eqref{T(v,1/v) as a double integral}
is slightly more involved.

    To prove Eq. \eqref{T(v,1/v) as a double integral}, first let $w=\frac{1}{v}$
in Eq. \eqref{T(v,w) as a double integral}.  This produces
\begin{equation}\label{T(v,1/v) proof part 1}
\begin{split}
\T\left(v,\frac{1}{v}\right)=&\tan^{-1}(v)\int_{0}^{1/v}\frac{\tan^{-1}(u)}{u}\d
u\\
&-\int_{0}^{\tan^{-1}(v)}\int_{0}^{\frac{\tan(\theta)}{v^2}}\frac{\tan^{-1}(z)}{z}\d
z\d\theta.
\end{split}
\end{equation}
Letting $v\rightarrow 1/v$ in Eq. \eqref{T(v,1/v) proof part 1}
gives
\begin{align*}
\begin{split}
\T\left(\frac{1}{v},v\right)=&\tan^{-1}\left(\frac{1}{v}\right)\int_{0}^{v}\frac{\tan^{-1}(u)}{u}\d
u-\int_{0}^{\tan^{-1}\left(\frac{1}{v}\right)}\int_{0}^{v^2\tan(\theta)}\frac{\tan^{-1}(z)}{z}\d
z\d\theta
\end{split}\\
\begin{split}
=&\left(\frac{\pi}{2}-\tan^{-1}(v)\right)\int_{0}^{v}\frac{\tan^{-1}(u)}{u}\d
u-\int_{\tan^{-1}(v)}^{\pi/2}\int_{0}^{\frac{v^2}{\tan(\theta)}}\frac{\tan^{-1}(z)}{z}\d
z\d\theta
\end{split}
\end{align*}
Now apply Eq. \eqref{arctangent integral functional equation}
twice, which transforms this last identity to
\begin{equation}\label{T(v,1/v) proof part 2}
\begin{split}
\T\left(\frac{1}{v},v\right)=&\left(\frac{\pi}{2}-\tan^{-1}(v)\right)
\left(\int_{0}^{\frac{1}{v}}\frac{\tan^{-1}(u)}{u}\d
u+\frac{\pi}{2}\log(v)\right)\\
&-\int_{\tan^{-1}(v)}^{\pi/2}\left(\int_{0}^{\frac{\tan(\theta)}{v^2}}\frac{\tan^{-1}(z)}{z}\d
z-\frac{\pi}{2}\log\left(\frac{1}{v^2}\tan(\theta)\right)\right)\d\theta
\end{split}
\end{equation}
To complete the proof, simply add equations \eqref{T(v,1/v) proof
part 1} and \eqref{T(v,1/v) proof part 2} together, and then
simplify the resulting sum. $\blacksquare$
\end{proof}

\begin{theorem}\label{T(v,1/v) mahler theorem}\allowdisplaybreaks{ Suppose that $v>0$, then the following Mahler
measures hold: {\allowdisplaybreaks
\begin{align}
\begin{split}
\m&\left(1-v^4\left(\frac{1-x}{1+x}\right)^2+\left(y+v^2\left(\frac{1-x}{1+x}\right)\right)^2
z\right)\\
&\qquad=\frac{4}{\pi}\int_{0}^{v}\frac{\tan^{-1}(u)}{u}\d
u-\frac{8}{\pi^2}\T\left(v,\frac{1}{v}\right)
+\frac{1}{2}\m\left(1-v^4\left(\frac{1-x}{1+x}\right)^2\right),
\end{split}\label{T(v,1/v) to mahler measure 1}\\
\begin{split}
\m&\left(1-v^4\left(\frac{1-x}{1+x}\right)^2+v^2\left(\frac{1-x}{1+x}\right)\left(\frac{1-y}{1+y}\right)
\left(z-z^{-1}\right)\right)\\
&\qquad=\frac{8}{\pi}\int_{0}^{v}\frac{\tan^{-1}(u)}{u}\d u
-\frac{16}{\pi^2}\T\left(v,\frac{1}{v}\right),
\end{split}\label{T(v,1/v) to mahler measure 2}\\
\begin{split}
\m&\left(\left(y-y^{-1}\right)+v^2\left(\frac{1-x}{1+x}\right)\left(z-z^{-1}\right)\right)\\
&\qquad=\frac{4}{\pi}\int_{0}^{v}\frac{\tan^{-1}(u)}{u}\d
u-\frac{8}{\pi^2}\T\left(v,\frac{1}{v}\right),
\end{split}\label{T(v,1/v) to mahler measure 3}
\end{align}
\begin{equation}\label{T(v,1/v) to mahler measure 4}
\begin{split}
\m&\left(\begin{split}&\left(4(1+y)^2-\left(z+z^{-1}\right)^2\right)\left(1-v^4\left(\frac{1-x}{1+x}\right)^2\right)^2\\
    &+\left(z-z^{-1}\right)^2(1+y)^2\left(1+v^4\left(\frac{1-x}{1+x}\right)^2\right)^2\end{split}\right)\\
&=\frac{8}{\pi}\int_{0}^{v}\frac{\tan^{-1}(u)}{u}\d
u-\frac{16}{\pi^2}\T\left(v,\frac{1}{v}\right)+\frac{4}{\pi}\int_{0}^{\pi/2}\log\left(1+v^2\tan(\theta)\right)\d\theta\\
&\quad+\log(2)
\end{split}
\end{equation}
}}
\end{theorem}
\begin{proof}Each of these results follows, in order, from substituting
Eq. \eqref{arctangent integral to mahler measure}, Eq.
\eqref{arctan extra mahler 1}, Eq. \eqref{arctan extra mahler 2},
and Eq. \eqref{arctan extra mahler 3}, into Eq. \eqref{T(v,1/v) as
a double integral}.$\blacksquare$
\end{proof}
\begin{corollary}
If we let $v=1$ in Eq. \eqref{T(v,1/v) to mahler measure 1}, we
can recover one of Lal\'in's formulas \cite{La1}:
\begin{align}
\m&\left((1+z)(1+y)+(1-z)(x-y)\right)=\frac{7}{2\pi^2}\zeta(3)+\frac{\log(2)}{2}.\label{lalin
formula} \end{align}
Letting $v=1$ in Eq. \eqref{T(v,1/v) to
mahler measure 2}, Eq. \eqref{T(v,1/v) to mahler measure 3}, and
Eq. \eqref{T(v,1/v) to mahler measure 4}, yields in order:
\begin{align}
\m&\left(4(1+y)+\left(1-y\right)\left(x-x^{-1}\right)\left(z-z^{-1}\right)\right)=\frac{14}{\pi^2}\zeta(3)\label{T(1,1) mahler 2}\\
\m&\left((1+x)\left(y-y^{-1}\right)+\left(1-x\right)\left(z-z^{-1}\right)\right)=\frac{7}{\pi^2}\zeta(3)\label{T(1,1)
mahler 3}\\
\begin{split}\m&\left(16(1+y)^2-4\left(z+z^{-1}\right)^2+(1+y)^2\left(z-z^{-1}\right)^2\left(x+x^{-1}\right)^2\right)\\
&\qquad=\frac{14}{\pi^2}\zeta(3)+\frac{4}{\pi}G\end{split}\label{T(1,1)
mahler 4}
\end{align}
\end{corollary}
\begin{proof}
To prove Eq. \eqref{lalin formula}, let $v=1$ in Eq.
\eqref{T(v,1/v) to mahler measure 1}.  From Eq. \eqref{T(v,1/v)
closed form} we know that
$\T(1,1)=\frac{\pi}{2}G-\frac{7}{8}\zeta(3)$, hence
\begin{align*}
\frac{7}{\pi^2}\zeta(3)+\log(2)&=\m\left(1-\left(\frac{1-x}{1+x}\right)^2+\left(y+\frac{1-x}{1+x}\right)^2
z\right)\\
&=\m\left(4x+\left((1+x)y+(1-x)\right)^2 z\right).
\end{align*}
Now let $(x,y,z)\rightarrow \left(x,\frac{y}{z},-x z^2\right)$ to
obtain
\begin{align*}
\frac{7}{\pi^2}\zeta(3)+\log(2)&=\m\left(4x-\left((1+x)y+(1-x)z\right)^2 x\right)\\
                               &=\m\left(4-\left((1+x)y+(1-x)z\right)^2\right)\\
                               &=2\m\left(2+(1+x)y+(1-x)z\right).
\end{align*}
With the final change of variables
$(x,y,z)\rightarrow\left(z,\frac{1}{y z},\frac{x}{yz}\right)$, we
have
\begin{align*}
\frac{7}{\pi^2}\zeta(3)+\log(2)&=2\m\left(2+\frac{(1+z)}{yz}+\frac{(1-z)x}{yz}\right)\\
                               &=2\m\left((1+z)(1+y)+(1-z)(x-y)\right),
\end{align*}
completing the proof of Eq. \eqref{lalin formula}.

    The proofs of Eq. \eqref{T(1,1) mahler 2} through Eq. \eqref{T(1,1) mahler 4} follow almost immediately from
our evaluation of $\T(1,1)$.  The proof Eq. \eqref{T(1,1) mahler
4} also requires the fairly easy fact that
$\int_{0}^{\pi/2}\log\left(1+\tan(\theta)\right)\d\theta=G+\frac{\pi}{4}\log(2)$
$\blacksquare$
\end{proof}
\section{Conclusion}
\label{conclusion} \init
    In principle, we should be able to apply the techniques
in this paper to prove formulas for infinitely many three-variable
Mahler measures.  The main difficulty, which is significant, lies
in the challenge of finding infinitely many Mahler measures for
the arctangent and arcsine integrals.  In Section \ref{mahler
measure integrals} we proved one such formula for the arcsine
integral, and four formulas for the arctangent integral.
%
%
\section{Acknowledgements}
I would like to thank my advisor, David Boyd, for bringing
Condon's paper to my attention.  He also made several useful
suggestions on evaluating the integrals in Theorem \ref{integral
table of jacobian integrals}.  I truly appreciate his support.

    Finally, I would like the thank the Referee for the helpful
suggestions.




\begin{thebibliography}{00}




\bibitem{Ba} N. Batir, Integral representations of some series
involving ${2k\choose k}^{-1}k^{-n}$ and some related series,
Applied Math. and Comp. \textbf{147} (2004), 645-667.

\bibitem{Be} B.C. Berndt, Ramanujan's Notebooks part IV,
Springer-Verlag, New York, 1994.

\bibitem{Bo0} D.W. Boyd, Speculations concerning the range of Mahler's measure,
Canad. Math. Bull. \textbf{24} (1981), 453-469.

\bibitem{Bo1} D.W. Boyd, Mahler's measure and special values of
L-functions, Experiment. Math. \textbf{7} (1998), 37-82.

\bibitem{BV} D.W. Boyd, F. Rodriguez Villegas, Mahler's measure and the dilogarithm (I),
Canad. J. Math. \textbf{54} (2002), 468-492.

\bibitem{Co} J. Condon, Calculation of the Mahler measure of a
three variable polynomial, (preprint, October 2003).

\bibitem{Gr} I.S. Gradshteyn, I.M. Ryzhik, Table
of Integrals, Series and Products, Academic Press, 1994.

\bibitem{La1} M. N. Lal\'{i}n, Some examples of Mahler measures as
multiple polylogarithms, J. Number Theory. \textbf{103} 2003,
85-108.

\bibitem{La2} M. N. Lal\'{i}n, Mahler measure of some $n$-variable polynomial families,
(preprint 2004, to appear in J. Number Theory).

\bibitem{La3} M. N. Lal\'{i}n, Some relations of Mahler measure with hyperbolic
volumes and special values of $\operatorname{L}$-functions,
(dissertation, 2005).

\bibitem{Le1} L. Lewin, Polylogarithms and Associated Functions, Elsevier North Holland, New York, 1981.

\bibitem{Ma} V. Maillot, G\'{e}om\'{e}trie d'Arakelov des vari\'{e}t\'{e}s toriques et
fibr\'{e}s en droites int\'{e}grables. M\'{e}m. Soc. Math. Fr.
(N.S) \textbf{80} (2000), 129pp.

\bibitem{Ra} S. Ramanujan, On the integral $\int_{0}^{x}\frac{\tan^{-1}(t)}{t}\d t$, J.
Ind. Math. Soc. \textbf{7} (1915), 93-96.

\bibitem{RV} F. Rodriguez Villegas, Modular Mahler measures I,
Topics in number theory (University Park, PA, 1997),
17--48, Math. Appl., 467, Kluwer Acad. Publ., Dordrecht, 1999.

\bibitem{Sm} C.J. Smyth, An explicit formula for the Mahler
measure of a family of 3-variable polynomials, J. Th. Nombres
Bordeaux, \textbf{14} (2002), 683-700.

\bibitem{Va} S. Vandervelde, A formula for the Mahler measure of
$a x y+bx+cy+d$, J. Number Theory, \textbf{100} (2003), 184-202.

\bibitem{Wa} G.N. Watson, A Treatise on the Theory
of Bessel Functions, Cambridge University Press, 1922.


\end{thebibliography}
\end{document}